\documentstyle[a4,amsfonts,12pt]{article}
\newcommand{\C}{{\bf C}}
\newcommand{\rbox}{$\:\:$ \raisebox{-1ex}{$\:\Box\:$}}
\newcommand{\OC}{\overline{{\bf C}}}
\newcommand{\OR}{\overline{{\bf R}}}
\newcommand{\R}{{\bf R}}
\newcommand{\He}{{\bf H}}

\newcommand{\D}{{\bf D}}
\newcommand{\cP}{{\cal P}}
\newcommand{\de}{\delta}

\newcommand{\mb}{\mbox}

\newcommand{\beq}{\begin{equation}}
\newcommand{\eeq}{\end{equation}}

\newcommand{\ve}{\varepsilon}

\newcommand{\ov}{\overline}
\newcommand{\al}{\alpha}
\newcommand{\be}{\beta}
\newcommand{\kap}{\mb{cap}}
\newcommand{\Om}{\Omega}
\newcommand{\Si}{\Sigma}

\newcommand{\z}{\zeta}

\newcommand{\ga}{\gamma}

\newcommand{\Ga}{\Gamma}
\newcommand{\si}{\sigma}
\newtheorem{th}{Theorem}

\newtheorem{lem}{Lemma}

\newcommand{\ueberschrift}{\bigskip\goodbreak\noindent\bigskip}
\newcounter{theabsatz}
\newcommand{\absatz}[1]{\stepcounter{theabsatz} \ueberschrift
                           {\large \bf \arabic{theabsatz}. {#1}} \setcounter{equation}{0}}

\begin{document}

\begin{center}
{\large \bf Positive Harmonic Functions on Denjoy Domains in the Complex Plane}\\[2ex]

by\\[2ex]

 Vladimir V.
Andrievskii\\[3ex]
{\it
Kent State University\\
Kent, OH, U.S.A.}
\end{center}


\footnotetext{


This research is supported in part by NSF grant DMS-0554344

 }

\vspace*{5mm}

\absatz{Introduction and Statement of Main Results}

 We use the following standard terminology.
  We denote by  Denjoy domain an
 open subset $\Om$ of the  complex plane $\C$
  whose complement
 $E:=\OC\setminus\Om$, where  $\OC:=\C\cup\{\infty\}$, is a subset
 of $\OR:=\R\cup\{\infty\}$, where $\R$ is the real axis (see \cite{garjon}).
  Throughout the paper we rely on the following assumption.
  Each point of $E$ (including the point at infinity)
 is regular for the Dirichlet problem in
 $\Om$.
  Denote by $\cP_\infty=\cP_\infty(\Om)$ the
 cone of positive harmonic functions on $\Om$ which have vanishing
 boundary values at every point of $E\setminus\{\infty\}$. Independently, Ancona
 \cite{anc} and Benedicks \cite{ben}  showed that
 either all functions in $\cP_\infty$ are proportional or
 $\cP_\infty$ is generated by two linearly independent (minimal)
 harmonic functions; that is, either dim $\cP_\infty=1$ or
 dim $\cP_\infty=2$ respectively. In other words, it means
 that the Martin boundary of $\Om$ has either one or two
 ``infinite" points.

 The results in \cite{anc} and \cite{ben} are proved for positive
 harmonic functions in domain $\Om\subset\R^n,n\ge2$.
 In this paper we focus on
  the case $n=2$ due to its  extreme importance in the theory
 of entire functions,
 where positive harmonic functions and
 subharmonic  functions in $\C$ which are non-positive on a subset of the real
 line were the subject of research significantly earlier.

  Bernstein \cite{ber23} showed that if an entire
 function $f$ satisfies
 \beq\label{1.1}
 \limsup_{|z|\to\infty}\frac{\log|f(z)|}{|z|}=:\sigma_f<\infty
 \eeq
 and
 $$
 |f(x)|\le1\quad (x\in\R),
 $$
 then $|f'(x)|\le\sigma_f$ for any $x\in\R.$

 In some extensions of the Bernstein theorem (see, for example,  \cite{pri16},
 \cite{akh46}, \cite{lev50}, \cite{sch}, \cite{akhlev}, \cite{lev71}, \cite{lev89}
 and  \cite{levlogsod}) the entire function $f$
 satisfies (\ref{1.1}) as well as a new condition
 $$
 |f(x)|\le1\quad (x\in E),
 $$
 where $E\subset\R$ conforms to certain metric properties. Then the
 authors derived estimates on the growth of $f$ in $\C$ of the form
 $$|f(z)|\le (H_E(z))^{\si_f}\quad (z\in\C),
 $$
 where $H_E(z)$ is a ``universal function" which does not
 depend on $f$.

 We say that a subharmonic function $u$ in $\C$ has finite degree
 $0<\sigma<\infty$, if
  $$
  \limsup_{|z|\to\infty}\frac{u(z)}{|z|}=\sigma.
  $$
  We denote by $K_\sigma(E)$ the class of subharmonic in $\C$
  functions of finite degree no larger than $\si$ and non-positive
  on $E$. Let
  $$v(z)=v(z,K_\si(E)):=\sup\{u(z):\, u\in K_\si(E)\}\quad
  (z\in\C)
  $$
  be the subharmonic majorant of the class $K_\si(E)$. It is known
  that $v(z)$ is either finite everywhere on $\C$ or equal to
  $+\infty$ on $\C\setminus E$. The set $E$ is said to be of type
  $(\al)$ in the former case, and of type $(\beta)$ in the latter.

  {\bf Theorem A.} (\cite{lev71}, \cite[Theorem 3.3]{lev89}, \cite[Remark 1]{levlogsod}) {\it

  (a) Case ($\al$) holds if and only if }dim {\it $\cP_\infty=2$;

  (b) case ($\beta$) holds if and only if }dim {\it  $\cP_\infty=1$.}

 There is  a close connection between the dimension of
 $\cP_\infty$ and the behavior of the Green function $g_\Om(\cdot,z)$
 for $\Om$ with pole at $z\in\Om$ (see \cite{tsu}, \cite{ran} or \cite{saftot} for further details on
 logarithmic potential theory). Let $E^*:=\R\setminus E$.

 {\bf Theorem B.} (\cite[Theorem 3]{kar}, \cite[Section VIII A.2]{koo})
 {\it Let
 $$U(z)=U(z,E):=\int_{E^*}g_\Om(t,z)dt\quad (z\in\Om).
 $$
 Then}

 (a) dim $\cP_\infty=1$  {\it if and only if $U\equiv\infty$
 in}
 $\Om$;

 (b) dim $\cP_\infty=2$  {\it if and only if $U$ is finite
 everywhere on
 $\Om$.}

 The problem of finding a geometric description of  $E$ such that
 dim $\cP_\infty=2$ or, equivalently, of  $E$ with the finite
 subharmonic majorants of classes $K_\sigma(E)$ attracted
 attention of a number of researches (see \cite{akhlev}, \cite{ben}, \cite{kar},
 \cite{seg88}, \cite{gar},
 \cite{seg90} and  \cite{sod}).

 One of the basic results in this area is the following Benedicks'
 criterion. Let
  $$R(x,r):=\left\{z\in\C:\,|\Re z-x|<\frac{r}{2},|\Im z|<\frac{r}{2}
  \right\}\quad (x\in\R,r>0).
  $$
  For an arbitrary fixed $\al$ with $0<\al<1$ and every
  $x\in\R\setminus\{0\}$, let $\be_x(\cdot)=\be_x(\cdot,\al,E
  )$ be the solution of the Dirichlet problem on
  $R(x,\al|x|)\setminus E$ with boundary values $\be_x=1$ on
  $\partial R(x,\al|x|)$ and $\be_x=0$ on $E\cap R(x,\al|x|)$.

  {\bf Theorem C.} (\cite[Theorem 4]{ben}) {\it For every $\al$ with
  $0<\al<1$,

  (a) $\mb{\em{ dim }}\cP_\infty=1$ if and only if
   $$\int_{|x|\ge1}\frac{\be_x(x)dx}{|x|}=\infty;
   $$

   (b) $\mb{\em{ dim }}\cP_\infty=2$ if and only if}
   $$\int_{|x|\ge1}\frac{\be_x(x)dx}{|x|}<\infty.
   $$
 Theorem C indicates that the dimension of
 $\cP_\infty$ only
 depends  on the geometry of $E$ near infinity.

 Theorems \ref{th1} and \ref{th2} below
 provide a natural and intrinsic characterization
  of  $E$ with a
 given
 dim $\cP_\infty$  in
 terms of the logarithmic capacity cap$(S)$, $S\subset\C$, which
 appears
 most suitable for this theory. In these theorems we also
 connect the dimension of $\cP_\infty$ with continuous properties
 of the Green function $g_\Om$ in a neighborhood of infinity.

 \begin{th}\label{th1}
  The following conditions are equivalent:

 (i) There exist points $a_j,b_j\in E,\, -\infty<j<\infty$ such
 that
  \beq\label{1.41n}
  b_{j-1}\le a_j<b_j\le a_{j+1},\quad \lim_{j\to\pm\infty}a_j=\pm\infty,
  \eeq
   \beq\label{1.42n}
  \bigcup_{j=-\infty}^\infty(a_j,b_j)\supset E^*,
  \eeq
  \beq\label{1.4}
  \inf_{-\infty<j<\infty}\frac{{\em \kap}(E\cap[a_j,b_j])}{{\em \kap}([a_j,b_j])}>0,
  \eeq
  \beq\label{1.5}
  \sum_{j=-\infty}^\infty\left(\frac{b_j-a_j}{|a_j|+1}\right)^2<\infty;
  \eeq

  (ii) {\em \mb{dim} }$\cP_\infty=2$;

  (iii) $\limsup_{\Om\ni t\to\infty}g_\Om(t,z)|t|<\infty$ for any
  $z\in\Om$.
 \end{th}

 For particular results of this kind,
 see \cite[Theorem 2]{seg88}, \cite[Theorem 8]{seg90}, \cite[Theorem
 4]{levlogsod} and \cite[Theorem 1.11]{cartot}.

 Notice that
 if $(a,\infty)\subset E^*$ or $(-\infty,a)\subset E^*$ for some
 $a\in\R$, then, by Theorem C, dim $\cP_\infty=1$.

 \begin{th}\label{th2} Let $E\cap(a,\infty)\neq\emptyset$ and
 $E\cap(-\infty,-a)\neq\emptyset$ for any $a>0$. The following conditions are equivalent:

 (i) There exist points $\{a_j,b_j\}_{j=-N}^M$, where
 $M+N=\infty$, such that $a_j,b_j\in E$,
  $$b_{j-1}\le a_j<b_j\le a_{j+1},$$
  \beq\label{1.6}
  \sup_{j}\frac{{\em \kap}(E\cap[a_j,b_j])}{{\em \kap}([a_j,b_j])}<1,
  \eeq
  \beq\label{1.7}
  \sum_{j=-N}^M\left(\frac{b_j-a_j}{|a_j+b_j|+1}\right)^2=\infty;
  \eeq

  (ii) {\em \mb{dim} } $\cP_\infty=1$;

  (iii) $\limsup_{\Om\ni t\to\infty}g_\Om(t,z)|t|=\infty$ for some
  $z\in\Om$.
 \end{th}

 Next, we state some metric tests which  immediately follow
 from the theorems above and which use
 the one-dimensional Lebesgue measure (length) $|S|$ of a linear
 (Borel)
 set $S\subset\R$.

 {\bf Remark 1.} Since
 $$\kap([a_j,b_j])=\frac{b_j-a_j}{4}\quad\mb{and}\quad
 \kap(E\cap[a_j,b_j])\ge\frac{|E\cap[a_j,b_j]|}{4},$$
 the existence of
 points $a_j,b_j\in E,\, -\infty<j<\infty$ satisfying
 (\ref{1.41n}), (\ref{1.42n}), (\ref{1.5}) and the new inequality
  \beq\label{1.8}
  \inf_{-\infty<j<\infty}\frac{|E\cap[a_j,b_j]|}{b_j-a_j}>0,
  \eeq
  is sufficient for the validity
 of the parts (ii) and (iii) of Theorem
  \ref{th1} (cf.  \cite[Lemma 1]{sch}, \cite[Theorem 5]{ben} and
   \cite[Theorem 3.8]{lev89}).

   {\bf Remark 2.} Let $E^*=\cup_{j=-N}^M(c_j,d_j)$, where
   $M+N=\infty$,
   be such
 that
  $$d_{j-1}<c_j<d_j<c_{j+1}.$$
 In Theorem \ref{th2}, taking the system of points $\{a_j,b_j\}$
 to be the same as $\{c_j,d_j\}$ we derive that
   the condition
   \beq\label{1.10}
   \sum_{j=-N}^M\left(\frac{d_j-c_j}{|c_j+d_j|+1}\right)^2=\infty
   \eeq
   implies each of the parts (ii) and (iii) of Theorem \ref{th2}
   (cf. \cite[Theorem 5]{ben}, \cite[Theorem 4]{kar} and \cite[Theorem 3.7]{lev89}).

  {\bf Remark 3.} Let
  $$\theta_E(t):=|E^*\cap[-t,t]|\quad (t>0).$$
  The condition
  $$
  \int_1^\infty\frac{\theta_E^2(t)dt}{t^3}<\infty
  $$
  yields each of the parts (ii) and (iii) of Theorem
  \ref{th1} (cf.  \cite[Theorem 4]{seg88}, \cite[Theorem 1]{gar},
  \cite[Theorem 13 and Theorem 17]{seg90} and
   \cite[Theorem 2.2 and Corollary 4.1]{tot}).

   {\bf Remark 4.} Let $\theta (t),t\ge 1$ be any
   increasing function such that
   \beq\label{1.hu1}
   0<\theta(t)\le 2t\quad (t\ge1),
   \eeq
   \beq\label{1.hu2}
  \int_1^\infty\frac{\theta^2(t)dt}{t^3}=\infty.
  \eeq
  Then, there exists $E$ such that
  \beq\label{1.3n}
  \theta_E(t)\le\theta(t)\quad (t>2),
  \eeq
   and parts (ii) and (iii) of Theorem \ref{th2} hold
   (cf. \cite[Section 4.4]{seg90}, \cite[Theorem 2]{gar} and
    \cite[Corollary 4.2]{tot}).

 {\bf Remark 5.}
 The continuous properties of the Green function $g_\Omega$ at boundary points
 are  of independent interest in potential
 theory (see, for example, \cite{car}, \cite{maz}, \cite{biavol}, \cite{and},
 \cite{cartot}, \cite{tot} and references therein).
 Using conformal invariance of the Green function and the linear transformation
  $$w=\frac{d_0-c_0}{2z-c_0-d_0},$$
   where $(c_0,d_0)$ is any of the finite components of $E^*$,  we can rephrase the part
 (i)$\Leftrightarrow$(iii) of Theorem \ref{th1} in the
 following equivalent form (in this form we  prove it in
 Section 5).

  Let $F\subset \R$ be a regular compact
set with the complement $G:=\OC\setminus F$, and let
 $g_G(\cdot)=g_G(\cdot,\infty)$ be the Green function of $G$ with pole
 at infinity.
 We assume that $F\subset[-1,1]=:I $ and  $\pm1,0\in
F$. Let $F^*:=\R\setminus F$.
 The equivalence (i)$\Leftrightarrow$(iii) in Theorem
  \ref{th1} can be restated as follows: {\it The following conditions
  are equivalent:

 (i') There exist points $a_j,b_j\in F,\, -\infty<j<\infty$ such
 that}
  $$
  -1\le a_{-1}<b_{-1}\le a_{-2}<b_{-2}\le\ldots<0<
  \ldots\le a_1 <b_1\le a_0<b_0\le 1,$$
 $$
  \bigcup_{j=-\infty}^\infty(a_j,b_j)\supset F^*\cap I,$$
  \beq\label{1.13}
  \inf_{-\infty<j<\infty}\frac{\kap(F\cap[a_j,b_j])}{\kap([a_j,b_j])}>0,
  \eeq
  \beq\label{1.14}
  \sum_{j=-\infty}^\infty\left(\frac{b_j-a_j}{a_j}\right)^2<\infty;
  \eeq

  (iii')
  \beq\label{1.15}
  \limsup_{G\ni z\to 0}\frac{g_G(z)}{|z|}<\infty.
  \eeq

The monotonicity of the Green function yields
 $$g_G(z)\ge
g_{\OC\setminus I}(z)\quad (z\in \C\setminus I),
$$
 that is, if $F$
has the ``highest density" at $0$, then $g_G$ has the ``highest
smoothness" at the origin. In particular,
 $$g_G(iy)\ge
g_{\OC\setminus I}(iy)>\frac{y}{2}\quad (0<y<1),$$ i.e.,
 $$ \limsup_{G\ni z\to 0}\frac{g_G(z)}{|z|}\ge\frac{1}{2}>0.
 $$
 In this regard, Remark 5
 describes the metric properties of $F$ such that $g_G$
has the ``highest smoothness" at $0$ (see the recent remarkable
result by Carleson and Totik \cite[Theorem 1.11]{cartot} for
another description of sets $F$ whose Green's function possesses
the property (\ref{1.15})).

 The rest of the paper is organized as follows.
 In Section 2 we  compile certain basic facts of the method
 by Akhiezer and Levin (see
 \cite{akhlev} and \cite{lev89}), who connected the positive harmonic
 functions in $\cP_\infty$ with majorants of classes of
 subharmonic functions via special conformal mapping  of the
 upper half-plane  onto the upper half-plane with vertical slits (a comb domain).
 Since a significant number of the proofs in this paper depends on a
 new technique for estimation of the module of families of
 curves, Section 3 contains a brief summary  of some special
 paths families and their modules.  In Section 4
 we  continue to discuss
 the modules of paths families (mainly crosscuts of a domain
 separating  subsets and points of its closure).
 Section 5 is devoted to the proof of Theorem
 \ref{th1}. In Section 6 we give the proof of Theorem \ref{th2}.  Remarks 3
 and 4 are proved in Section 7.

\absatz{Majorants  in Classes of Subharmonic Functions}

 To
describe the properties of  subharmonic majorants of the class
$K_\si(E)$, Akhiezer and Levin \cite{akhlev} introduced a special
conformal mapping of the upper half-plane $\He:=\{z:\, \Im z>0\}$
onto $\He$ with vertical slits. Later, Levin \cite{lev89}
constructed the general theory of such mappings and used them to
solve several extremal problems in classes of subharmonic
functions.
 In this section we discuss some basic results of this theory which we use in the
 proofs of Theorem \ref{th1} and Theorem \ref{th2}.

 Recall that according to Theorem C we can assume that $0$ belongs
 to $E$ with some interval around the origin.
 Let $E^*=\R\setminus E=\cup_{k=0}^M(c_k,d_k)$, where
 $c_k<d_{k}$, $(c_k,d_k)\cap(c_j,d_j)=\emptyset$ if $k\neq j$ and
 $M\ge 0$ can
 be either finite or infinite.  As before we assume that each point of $E$ is regular
 for the Dirichlet problem in $\Om$ and that
 $E\cap(a,\infty)\neq\emptyset$ and
 $E\cap(-\infty,-a)\neq\emptyset$ for any $a>0$.
 By
 \cite[Theorem 3.1]{lev89} there exists a conformal mapping
 $\phi(\cdot)=\phi(\cdot,E)$ of $\He$ onto
  $$\He_E=\He\setminus \bigcup_{k=0}^M(u_k,u_k+iv_k],
  $$
  where $u_k\in\R,u_k\neq u_j$ for $k\neq j$ and $v_k>0$, such
  that
  $$\phi(0):=\lim_{\He\ni z\to 0}\phi(z)=0,
  \quad \phi(\infty):=\lim_{\He\ni z\to \infty}\phi(z)=\infty,$$
  $$\phi([c_k,d_k]):=\bigcup_{c_k\le x\le d_k}
  \lim_{\He\ni z\to  x}\phi(z)=[u_k,u_k+iv_k],
  $$
   and any point of $[u_k,u_k+iv_k)$ has exactly
  two preimages. Besides,
  $$\phi(E):=\bigcup_{x\in E}\lim_{\He\ni z\to
  x}\phi(z)=\OR.
  $$
  Note that $\phi$ is defined up to the multiplication by a
  positive constant. As a shorthand, we sometimes use
 $\phi$ both for the conformal mapping of $\He$ and for
  its continuous extension to $\ov{\He}$.
  \begin{lem}\label{lem2.1}
  (Theorem A and \cite[Section $n^0 5$]{akhlev}). \em{ dim }$\cP_\infty=2$
  if and only if
 $$
  \limsup_{y\to +\infty}\frac{|\phi(iy)|}{y}>0.
  $$
  \end{lem}
  Along with the set $E$, we consider the family of sets
  $$E_R:=E\cup\{ x\in\R:\, |x|\ge R\}\quad (R>0).
  $$
  Let $\phi_R(\cdot):=\phi(\cdot,E_R)$ be the appropriate conformal
  mapping. Multiplying by a positive constant if necessary, we
  normalize $\phi_R$ such that
  \beq\label{2.2}
  \lim_{\He\ni z\to\infty}\frac{\phi_R(z)}{z}=1.
  \eeq
 \begin{lem}
 \label{lem2.2}
 (Theorem A, \cite[Section $n^0 3$]{akhlev}, and
 \cite[Theorem  3]{levlogsod}). \em{ dim }$\cP_\infty=2$ if and only if $\Im
 \phi_R(z)=O(1)$ as $R\to\infty$ at some $z\in\He$.
 \end{lem}

 In the case of a regular (for the Dirichlet problem) compact set $F\subset\R$ we use an
 analogue of a conformal mapping $\phi$ to describe the Green
 function $g_G(\cdot)=g_G(\cdot,\infty)$ for $G=\OC\setminus F$ (with pole at $\infty$)
 and the capacity of $F$ (see \cite{and04} for details).
 Applying linear transformation if necessary we can always  assume
 that $F\subset[-1,1]=:I$ and
 $\pm 1\in F$. Consider the nontrivial case where $F\neq I$. The open (with respect to
$\R$) set $F^*:=\R\setminus F$ consists of either a finite number
$N>2$ or infinite number $N=\infty$ of open disjoint intervals,
i.e.,
 $$F^*=\bigcup_{j=1}^N(c_j,d_j),$$
  where
$(c_1,d_1):=(-\infty,-1)$ and $ (c_2,d_2):=(1,\infty).$
 Let
$\mu_F$ be the equilibrium measure for the set $F$.
 Consider the function $$
\psi(z)=\psi(z,F):=\pi+i\left(\int_F\log(z-\z)\,
d\mu_F(\z)-\log\kap(F)\right) \quad (z\in \He).$$
 It is analytic and univalent in
$\He$ and it maps $\He$
 onto  a vertical half-strip
with $N-2$  slits parallel to the imaginary axis, i.e., the domain
 $$\Sigma_F=\{ w:\, 0<\Re w<\pi,\Im w>0\}\setminus
\bigcup_{j=3}^{N}[\tilde{u}_j,\tilde{u}_j+i\tilde{v}_j],$$
 where
$0<\tilde{u}_j<\pi$ and $\tilde{v}_j>0.$

 The continuous extension of $\psi$ to $\ov{\He}$ satisfies the following
boundary correspondence
$$\psi(\infty)=\infty,\, \psi((-\infty,-1])=\{w:\, \Re w=0,\Im w\ge
0\},$$
 $$ \psi([1,\infty))=\{w:\, \Re w=\pi,\Im w\ge 0\},\quad
  \psi(F)=[0,\pi],$$
 $$\psi([c_j,d_{j}])=[\tilde{u}_j,\tilde{u}_j+i\tilde{v}_j]\quad (j=3,\ldots,N).$$
 Note
that in the last relation above, each point of
$[\tilde{u}_j,\tilde{u}_j+i\tilde{v}_j)$ has two preimages on
$[c_j,d_{j}]$.

 The advantage of using $\psi$ lies in the fact  that
  $$
g_G(z)=\Im\{\psi(z)\}\quad (z\in \ov{\He}).
 $$

 Consider the function $\Psi(z):=e^{i(\pi-\psi(z))}, \, z\in\ov{\He}$.
 Using the reflection principle we  extend
$\Psi$ to a function analytic in $\OC\setminus I$ according to the
formula
 $$\Psi(z):=\ov{\Psi(\ov{z})} \quad (z\in\C\setminus
\ov{\He}).
 $$
  Let ${\D}:=\{w:\, |w|<1\}$ be the unit disk. The function
 $\Psi$ is analytic and univalent and it maps $\OC\setminus I$ onto a starlike
(with respect to $\infty$) domain $\OC\setminus K_F$ with the
following properties: it is symmetric with respect to the real
line $\R$ and it coincides with the exterior of the unit disk with
$2N-4$ slits, i.e.,
 $$\ov{\He}\cap
K_F=(\ov{\He}\cap\ov{\D})\cup(\bigcup_{j=3}^{N}
[e^{i\theta_j},r_je^{i\theta_j}]),
 $$
  where $r_j:=e^{\tilde{v}_j}>1$ and
 $0<\theta_j:=\pi-\tilde{u}_j<\pi$.

  There is a close connection
between the capacities of the compact sets $K_F$ and $F$, namely
 \beq\label{2.3}
  \kap(K_F)=\frac{1}{2\,\kap(F)}=\frac{\kap(I)}{\kap(F)}\, .
  \eeq

\absatz{Modules of Path Families}

The main idea of our approach is to estimate the (equal) modules
of a family $\Ga$ of some paths in $\He$ and the family
$\phi(\Ga)$ of paths in $\He_E$ or the family $\psi(\Ga)$ in
$\Si_F$ in various ways. We briefly recall the notion of the
module of families of curves, generalized in an obvious way to
path families, see \cite{ahl}, \cite{lehvir},  \cite{pom} and
\cite{garmar} for details.

As usual, a Jordan curve is a continuous image of a closed
interval without intersections (except possibly endpoints). By a
curve we understand a locally rectifiable Jordan curve without
endpoints. We define a path to be the union of finitely many
mutually disjoint curves. For any path $\ga$, denote by $\ov{\ga}$
the closure of $\ga$ in $\C$, i.e., the union of $\ga$ and the
endpoints of curves composing $\ga$.

A Borel measurable function $\rho\ge 0$ on $\C$ is a metric if
 $$0<A(\rho):=\int\limits_{\C}\rho^2(z)\, dm(z)<\infty,$$
  where
$dm(z)$ stands for the $2$-dimensional Lebesgue measure (area) on
$\C$.

For a path family $\Ga=\{\ga\}$  let
$$L_\rho(\Ga):=\inf_{\ga\in\Ga}\int\limits_\ga\rho(z)\, |dz|$$ (if
the latter integral does not exist for some $\ga\in\Ga$, then we
define it to be the infinity).

The quantity
 \beq\label{aaa}
m(\Ga):=\inf_{\rho}\frac{A(\rho)}{L_\rho^2(\Ga)}\, ,\eeq
 where the
infimum is taken with respect to all metrics $\rho$, is called
 the module of the family $\Ga$.
  In the sequel we refer to the basic properties of the module, such as
conformal invariance, comparison principle, composition laws, etc.
 (see \cite{ahl}, \cite{lehvir}, \cite{pom} and \cite{garmar}). We
will use these properties and the definition (\ref{aaa}) without
further citations.

Special families of separating paths play a rather useful role.
Let $D\subset\OC$ be a domain. We say that a path $\ga\in D$
separates sets $A\subset\ov{D}$ and $B\subset\ov{D}$ if any Jordan
curve $J\subset D$ joining $A$ with $B$ has nonempty intersection
with $\ga$.

We use $\Ga,\Ga_1,\ldots$ to denote path families. We may use the
same symbol for different families if it does not lead to
confusion.

The  examples below state  some  well-known facts concerning
special families of curves.

For $w_0\in\R$ and $0<r_1<r_2,$ let
$$\Ga_1=\Ga_1(w_0,r_1,r_2):=\{\ga_r=\{w_0+re^{i\theta}:0<\theta<\pi\}:r_1<r<r_2\},$$
 $$\Ga_2=\Ga_2(w_0,r_1,r_2):=\{\ga_\theta=\{w_0+re^{i\theta}:r_1<r<r_2\}:0<\theta<\pi\}.$$
 Then
 \beq\label{3.nn4}
 m(\Ga_1)=\frac{1}{\pi}\log\frac{r_2}{r_1},
 \eeq
  \beq\label{3.nn5}
 m(\Ga_2)=\frac{\pi}{\log\frac{r_2}{r_1}}.
 \eeq
 Further, for $a<b$ and $c>0$, let $\Ga_3=\Ga_3(a,b,c)$ be the family of all
 half-ellipses in $\He$ with the foci at $a$ and $b$ which have
 nonempty intersection with the interval
 $\left(\frac{a+b}{2},\frac{a+b}{2}+ic\right)$.
 Applying the
 transformation
 $$w\to\frac{b-a}{4}\left(w+\frac{1}{w}\right)+\frac{a+b}{2}$$
 and (\ref{3.nn4}) we obtain
 \begin{eqnarray}
 m(\Ga_3)&=&\frac{1}{\pi}\log
 \left(\frac{2c}{b-a}+\sqrt{1+\frac{4c^2}{(b-a)^2}}\right)\nonumber\\
 \label{3.nn3}
 &>& \frac{1}{2\pi}\log
 \left(1+\frac{4c^2}{(b-a)^2}\right).
 \end{eqnarray}
 Let $\Ga_4=\Ga_4(r,R),\, 0<r<R$ be the family of all closed
curves in $\He$ which separate points $ir$ and $iR$ from $\R$.
Applying the linear transformation
 $$z\to\frac{z-ir}{z+ir},$$
 we see that
 $$m(\Ga_4)=\frac{1}{2\pi}\mu\left(\frac{R-r}{R+r}\right)=
 \frac{\pi}{4\mu\left(\frac{r}{R}\right)},$$
 where $\mu(t),\, 0<t<1$ is the module of Gr\"otzsch's
 extremal domain $\D\setminus[0,t]$ (see \cite[pp. 53,
 61]{lehvir}). Since
 \beq\label{3.k1}
 \mu(t)<\log\frac{4}{t}
 \eeq
 (see \cite[p. 61]{lehvir}), it follows that
 \beq\label{3.1}
 m(\Ga_4)\ge\frac{\pi}{4\log\frac{4R}{r}}.
 \eeq
 Further, denote by $\Ga_5=\Ga_5(r,R),\, 0<r<R$, the family of
 all crosscuts of $\He$, i.e., curves in $\He$ with endpoints on $\R$,
  which separate $0$ and $ir$ from $i R$ and
 $\infty$. By  symmetry
 $$2m(\Ga_4)2m(\Ga_5)=1.$$
 Therefore, (\ref{3.1}) implies that
 $$
 m(\Ga_5)=\frac{1}{4m(\Ga_4)}\le\frac{1}{\pi}\log\frac{4R}{r}.
 $$
For $a<b<c$, let $\Ga_6=\Ga_6(a,b,c)$ denote the family of all
crosscuts of $\He$ joining the boundary intervals $(a,b)$ and
$(c,\infty)$. The  appropriate result for Teichm\"uler's extremal
problem (see \cite[p. 55]{lehvir}) and (\ref{3.k1}) yield
 \beq\label{3.nn1}
 m(\Ga_6)\le\frac{1}{\pi}\log\frac{16(c-a)}{c-b}.
 \eeq
 For $a>0$, let $\Ga_7=\Ga_7(a)$ be the family of all curves
 joining the boundary intervals $[0,ia]$ and $\{z:\, \Re z=\pi,\Im
 z>0\}$ in the half-strip $\{z:\, 0<\Re z<\pi,\Im z>0\}.$
 By \cite[(4.8)-(4.9)]{and} we have
  $$m(\Ga_7)\ge \left\{\begin{array}{ll}
\frac{1}{2},&\mb{ if }a\ge\frac{\log 2}{2},\\[2ex]
 \frac{\pi}{2\log \frac{12}{a}}&\mb{ if } a\le\frac{\log 2}{2}.
\end{array}\right.
 $$
 For $a\in\R$ and $0<b<c$, let $\Ga_8=\Ga_8(a,b,c)$ be the family
 of all curves $\ga\subset\Sigma=\Sigma(a,b,c):=\{z:\, a-\pi<\Re z<a+\pi,\Im z>b\}$
  joining $\{z:\, \Re z=a-\pi,\Im z>b\}$ and $\{z:\, \Re z=a+\pi,\Im z>b\}$
  in $\Sigma$  such that $\ga\cap \{z:\, \Re z=a,\Im z>b\}$
  consists of exactly one point which belongs to the interval
  $(a+ib,a+ic)$. The result of the previous example implies that
   \beq\label{3.51}
   m(\Ga_8)\ge \left\{\begin{array}{ll}
\frac{1}{4},&\mb{ if }c-b\ge\frac{\log 2}{2},\\[2ex]
 \frac{\pi}{4\log \frac{12}{c-b}}&\mb{ if } c-b\le\frac{\log 2}{2}.
\end{array}\right.
 \eeq
 For $0<a<b$ and $c>0$, denote by $\Ga_9=\Ga_9(a,b,c)$ the family
 of all crosscuts of
 $$G=G(a,b,c)=\He\setminus R_1,$$
 where $R_1=R_1(a,b,c):=\{z:\, a\le\Re z\le b,0<\Im z\le c\}$, joining
 boundary intervals $(0,a)$ and $(b,\infty)$. We claim that
 \beq\label{3.nn2}
 m(\Ga_9)\ge\frac{1}{2\pi}\log\frac{b^2}{c^2+(b-a)^2}.
 \eeq
 Indeed, in the nontrivial case where  $c^2+(b-a)^2<b^2$ we apply (\ref{3.nn4})
 to have
 $$m(\Ga_9)\ge m(\Ga_1(b,\sqrt{c^2+(b-a)^2},b)= \frac{1}{2\pi}\log\frac{b^2}{c^2+(b-a)^2}.
 $$
 Next, we discuss the following direct
 consequence  of Pfl\"uger's theorem (see \cite[p. 212]{pom}).
  Let $0<a<b$ be such that
 $\al:=\log\frac{b}{a}<\frac{\pi}{2},$ and let $S\subset[a,b]$ be any
   set consisting of a finite number of closed intervals. Denote
   by $\Ga_{10}=\Ga_{10}(a,b,A)$ the family of
 all paths in
 $$R_2=R_2(a,b):=\{re^{i\theta}:\, a<r<b,0<\theta<\al\}$$
 which separate $S$ from the interval $(ae^{i\al},be^{i\al}):=
 \{re^{i\theta}:\, a<r<b,\theta=\al\}$.
 We claim that
 \beq\label{3.3}
 m(\Ga_{10})\le \frac{2}{\pi}\log\frac{C(b-a)}{\kap(S)},
 \quad C=\frac{e^\pi+1}{2}.
 \eeq
 To see this,
 consider the function
 $$f(z)=\exp\left\{i\pi\frac{\log\frac{z}{a}}{\log\frac{b}{a}}\right\}
 $$
 which maps $R_2$ conformally and univalently onto the half-ring
 $$R_3=\{\z\in\He:\, e^{-\pi}<|\z|<1\}.$$
 Since
 $$|f(x_2)-f(x_1)|\ge\frac{2|x_2-x_1|}{\al b}\quad (x_1,x_2\in S),$$
 for $B:=f(S)$, we have
 \beq\label{3.4}
 \kap(B)\ge \frac{2}{\al b}\kap(S).
 \eeq
 Let $B':=\{z\in\C\setminus\He:\,\ov{z}\in B\}.$ Pfl\"uger's
 inequality
  (see \cite[p. 212]{pom}) yields
  \beq\label{3.k2}
  \kap(B)\le\kap(B\cup B')\le
  (e^{\pi/2}+e^{-\pi/2})\exp\left\{-\frac{\pi}{ m(\Ga_{11})}\right\},
  \eeq
  where $\Ga_{11}$ is the family of all curves in the ring $\{z:\,
  e^{-\pi}<|z|<1\}$  joining $B\cup B'$ with $\{z:\,
  |z|=e^{-\pi}\}$.

  By symmetry  and \cite[Theorem IV.4.2]{gar}
  $$m(\Ga_{11})=\frac{2}{m(\Ga_{10})}.$$
   Therefore, (\ref{3.4}) and (\ref{3.k2}) imply (\ref{3.3}).

 Next, let
 $$-R\le\al_{-M}<u_{-M}<\be_{-M}\le \ldots<\be_{-2}
 \le \al_{-1}<u_{-1}<\be_{-1}\le-r<0$$
 $$<r\le\al_0<u_0<\be_0\le\al_1<
 \ldots\le\al_N<u_N<\be_N\le R,
 $$
 where $M>0$ and $N\ge0$ are finite integers, and let  real numbers $u_j$ and
 $v_j$,  $-M\le j\le N$ satisfy
 $$0\le
 v_j\le \min\{|u_j|,\be_j-\al_j\}
 .$$
 Denote by $\Ga_{12}=\Ga_{12}(r,R,\{u_j\},\{v_j\})$ the family of all curves in
 $$D=D(r,R,\{u_j\},\{v_j\}):=\{w\in\He:\, r<|w|<2R\}\setminus
 \bigcup_{j=-M}^{N}
 [u_j,u_j+iv_j]$$
 which join circular arcs $\{z\in\He:\,|z|=r\}$ and $\{z\in\He:\,|z|=2R\}$.
 \begin{lem}\label{lema}
 The inequality
 \beq\label{3.5}
 m(\Ga_{12})\le\frac{\pi\log\frac{2R}{r}-C\sum_{j=-M}^{N}
 \left(\frac{v_j}{u_j}\right)^2}{
 \left(\log\frac{2R}{r}\right)^2}
 \eeq
 holds with the  constant $C=10^{-5}$.
 \end{lem}
  {\it Proof.}
 Let
$\ve:=10^{-1}$, and for $-M\le j\le N$, let
 $$B_j:=
 \left\{\begin{array}{ll}
 \{ w\in D:\, |w-u_j|\le \ve^2 v_j,\, \Re
w\le u_j\},& \mb{ if }u_j\ge\frac{1}{2}(\be_j +\al_j),\\[2ex]
 \{ w\in D:\, |w-u_j|\le \ve^2 v_j,\, \Re
w\ge u_j\},& \mb{ if }u_j<\frac{1}{2}(\be_j
+\al_j).\end{array}\right.
$$
 Consider the metric
 $$ \rho(w)=
\left\{\begin{array}{ll} |w|^{-1}\, ,& \mb{ if }w\in
 D\setminus\cup_{j=-M}^{N}B_j,\\[2ex] 0,&\mb{
 elsewhere}.\end{array}\right.$$
 We proceed to show that
 \beq\label{5.9n}
 \int\limits_\ga\rho(w)|dw|\ge\log\frac{2R}{r}\quad
 (\ga\in\Ga_{12} ) .
 \eeq
 Indeed,
 if $\ga\cap
 (\cup_{j=-M}^{N}B_j)=\emptyset$, then
 $$
 \int\limits_\ga\rho(w)|dw|\ge\left|\int\limits_{\ga}
\frac{dw}{w}\right|\ge\log\frac{2R}{r}\, .
$$
 Moreover,  (\ref{5.9n}) remains  valid even if $\ga\cap(\cup_{j=-M
 }^{N}B_j)\neq\emptyset.$ We prove this as follows. Let $j$ be such that
 $\ga\cap B_j\neq\emptyset$.
 Denote by $\ga_j\subset\ga$ the largest curve which has nonempty
 intersection with $B_j$ and whose endpoints $t_j$ and $\tau_j$
 are such that $\ve t_j,\ve\tau_j\in \{w\in\partial B_j:\,
 |w-u_j|=\ve^2
 v_j\}$. Let $\ga'_j\subset D\cap\{w:\,|w-u_j|=\ve v_j\}$ be the
 circular arc with the same endpoints $t_j$ and $\tau_j$.
 Since $\ga_j$ includes two curves joining $\{w:\,
 |w-u_j|=\ve^2 v_j\}$ and $\{w:\,
 |w-u_j|=\ve v_j\}$, we conclude that
 $$\int_{\ga_j}\rho(w)|dw|\ge\frac{2\ve(1-\ve)v_j}{|u_j|+\ve v_j}.$$
 Meanwhile, since $\ga'_j$ is a subarc of a circular part
 of  $\partial B_j$, we obtain
 $$\int_{\ga'_j}\rho(w)|dw|\le\frac{\pi}{2}\frac{\ve v_j}{|u_j|-\ve v_j}.$$
 Hence, the choice of sufficiently small $\ve=10^{-1}$ guaranties
 that
 $$\int_{\ga_j}\rho(w)|dw|\ge \int_{\ga'_j}\rho(w)|dw|.$$
 Therefore,
 $$\int_{\ga}\rho(w)|dw|\ge \int_{\ga'}\rho(w)|dw|\ge\log\frac{2R}{r},$$
 where
 $$\ga'=\left(\ga\setminus\bigcup_{\ga\cap
 B_j\neq\emptyset}\ga_j\right)\cup \left(\bigcup_{\ga\cap
 B_j\neq\emptyset}\ga_j'\right)\subset D\setminus\bigcup_{j=-M}^NB_j.
 $$
 This completes the proof of (\ref{5.9n}).

 Since
  $$A(\rho)\le\pi\log\frac{2R}{r}-\frac{\pi}{4}\sum_{j=-M
 }^{N}
 \frac{\ve^4v_j^2}{(|u_j|+\ve^2v_j)^2}\, ,$$
 according to (\ref{5.9n}) we have (\ref{3.5}).

 \hfill\rbox

 Further, let for $T>1$,
 $$-T\le a_{-M}<b_{-M}\le a_{-M+1}<\ldots\le a_{-1}<b_{-1}\le
 -1,$$
 and
 $$1\le a_0<b_0\le a_1<\ldots\le a_N< b_N\le T,$$
 where $M>0$ and $N\ge 0$ are finite, be such that
  $$\max_{-M\le j\le N}\frac{b_j-a_j}{\min\{|a_j|,|b_j|\}}<\frac{\pi}{2}.$$
 Let
  $$S^*\subset\bigcup_{j=-M}^N(a_j,b_j)$$
  consist of a finite number of open (disjoint) intervals
  such that $S:=([-T,-1]\cup[1,T])\setminus
  S^*$ is a regular set satisfying
  \beq\label{3.11n}
  \min_{-M\le j\le
  N}\frac{\kap(S\cap[a_j,b_j])}{\kap([a_j,b_j])}\ge q>0.
  \eeq
  Denote by $\Ga_{13}=\Ga_{13}(T,S,\{(a_j,b_j)\})$ the family of all paths
  in
  $$Q=Q(T):=\{z\in\He:\, 1<|z|<T\}
  $$
  which separate either
  $S\cap[-T,-1]$ from $[1,T]$ or $S\cap[1,T]$ from $[-T,-1]$.
  \begin{lem}\label{lem3.1n}
  Under the  above assumptions  the inequality
 \beq\label{3.12n} m(\Ga_{13})\le\frac{\pi\log T+C\sum_{j=-M}^N
 \left(\frac{b_j-a_j}{\min\{|a_j|,|b_j|\}}\right)^2}{\left(\log T\right)^2}\, ,
 \eeq
  holds with the constant $C=\frac{2}{\pi}\log\frac{10^2}{q}$.
  \end{lem}
  {\it Proof.}
  There is no loss of generality in assuming that
  $S^*\cap(a_j,b_j)\neq\emptyset$ for any $-M\le j\le N$.
 Let
 $$R_j:=\left\{\begin{array}{ll}
\{z=re^{i\theta}:\, a_j\le r\le
b_j,0\le\theta\le\log\frac{b_j}{a_j}\},&\mb{ if }0\le j\le N,\\[2ex]
 \{z=re^{i\theta}:\, -b_j\le r\le
-a_j,\pi-\log\frac{a_j}{b_j}\le\theta\le\pi\},&\mb{ if }-M\le j\le
 -1.
\end{array}\right.
$$
 Denote by $\Ga_j$ the family of all paths in $R_j$ which separate
 $S\cap[a_j,b_j]$ and the boundary interval
 $$\left\{\begin{array}{ll}
\{z=re^{i\theta}:\, a_j\le r\le
b_j,\theta=\log\frac{b_j}{a_j}\},&\mb{ if }0\le j\le N,\\[2ex]
 \left\{z=re^{i\theta}:\, -b_j\le r\le
-a_j,\theta=\pi-\log\frac{a_j}{b_j}\right\},&\mb{ if }-M\le j\le
-1.
 \end{array}\right.
 $$
 According to (\ref{3.3}) and our assumption (\ref{3.11n})   for $j$ under
  consideration,  there exists
 the metric $\rho_j$ (with the support in $R_j$) such that
 $$\int_\ga \rho_j(z)|dz|\ge 1\quad (\ga\in\Ga_j)$$
 and
 \beq\label{3.v1}
 A(\rho_j)\le\frac{2}{\pi}\log\frac{(e^\pi+1)(b_j-a_j)}
 {2\kap(S\cap[a_j,b_j])} \le\frac{2}{\pi}\log\frac{10^2}{q}=:C.
 \eeq
 Consider the metrics
 $$\rho^*(z):=\left\{\begin{array}{ll}
|z|^{-1},&\mb{ if }z\in Q,\\[2ex]
 0,&\mb{ elsewhere},
\end{array}\right.$$
 $$\rho(z):=\max\left\{\rho^*(z),\sum_{j=-M }^N
  \rho_j(z)\left|\log
 \frac{a_j}{b_j}\right|\right\}.$$
 We claim that for any $\ga\in\Ga_{13}$,
 \beq\label{3.13n}
 \int_\ga\rho(z)|dz|\ge\log T.
 \eeq
 Below we confirm the validity of (\ref{3.13n}) in the case where
 $\ov{\ga}\cap[-T,-1]=\emptyset$
 (the proof of
 (\ref{3.13n}) in the case where $\ov{\ga}\cap[1,T]=\emptyset$
 follows along the same lines).
 Let for  $\ga\in\Ga_{13}$ and $1\le r_1\le r_2\le T$,
 $$\ga(r_1,r_2):=\{ z\in\ga:\, r_1\le |z|\le r_2\}.
 $$
 If $\ov{\ga(r_1,r_2)}\cap[r_1,r_2]=\emptyset$, then
 \beq\label{3.112n}
 \int_{\ga(r_1,r_2)}\rho(z)|dz|\ge
 \int_{\ga(r_1,r_2)}\rho^*(z)|dz|\ge \left|\int_{\tilde{\ga}(r_1,r_2)}\frac{dz}{z}
 \right|\ge\log\frac{r_2}{r_1}\, ,
 \eeq
 where $\tilde{\ga}(r_1,r_2)\subset \ga(r_1,r_2)$ is any curve
 joining $\{z:\, |z|=r_1\}$ with $\{z:\, |z|=r_2\}$.

 Furthermore, let $\ov{\ga(a_j,b_j)}\cap[a_j,b_j]\neq\emptyset$ for some $0\le
 j\le N$. If  there is no  curve
 $\tilde{\ga}(a_j,b_j)\subset \ga(a_j,b_j)$ joining the circles $\{
 z:\,|z|=a_j\}$ with $\{
 z:\,|z|=b_j\}$,
 then at least one of the following two cases holds.

 (a) If there exists $\ga_j'\subset\ga(a_j,b_j)$ such that
 $\ga_j'\in\Ga_j$, then
 \beq\label{3.113n}
 \int_{\ga(a_j,b_j)}\rho(z)|dz|\ge\left(\log\frac{b_j}{a_j}\right)\int_{\ga_j'}\rho_j(z)|dz|\ge
 \log\frac{b_j}{a_j}.
 \eeq

 (b) If there exists $\ga_j''\subset\ga(a_j,b_j)\cap R_j$ joining
 $\{z:\, \arg z=\log\frac{b_j}{a_j}\}$ with $\R$, then
 \beq\label{3.114n}
 \int_{\ga(a_j,b_j)}\rho(z)|dz|\ge\int_{\ga_j''}\rho^*(z)|dz|\ge
 \left|\int_{\ga_j''}\frac{dz}{z}\right|\ge
 \log\frac{b_j}{a_j}.
 \eeq
 Therefore, by (\ref{3.112n})-(\ref{3.114n})  we obtain
 \begin{eqnarray*}
 \int_\ga\rho(z)|dz|&\ge&\sum_{j=0}^N\int_{\ga(a_j,b_j)}\rho(z)|dz|+
 \sum_{j=0}^{N-1}\int_{\ga(b_j,a_{j+1})}\rho(z)|dz|\\
 &+& \int_{\ga(1,a_0)}\rho(z)|dz|
 +\int_{\ga(b_N,T)}\rho(z)|dz|\\
 &\ge&\sum_{j=0}^{N}\log\frac{b_j}{a_j}+
 \sum_{j=0}^{N-1}\log\frac{a_{j+1}}{b_j}+\log a_0+\log\frac{T}{b_N}=\log T,
 \end{eqnarray*}
 which proves (\ref{3.13n}).

 Since by (\ref{3.v1})
 \begin{eqnarray*}
 A(\rho)&\le& A(\rho^*)+\sum_{j=-M}^NA(\rho_j)\left(\log\frac{b_j}{a_j}\right)^2\\
 &\le&\pi\log T+C\sum_{j=-M}^N
 \left(\frac{b_j-a_j}{\min\{|a_j|,|b_j|\}}\right)^2,
 \end{eqnarray*}
 by virtue of  (\ref{3.13n})  and the definition of the module we have (\ref{3.12n}).

 \hfill\rbox

\absatz{Preliminary Constructions}

In this section we prove some auxiliary results.

Let $F\subset I,F\neq I,\pm 1\in F,F^*,\psi:\C\setminus
I\to\Sigma_F,$ and $\Psi:\OC\setminus I\to\OC\setminus K_F$ be
defined as in Section 2, and let $(c_j,d_j),3\le j\le N\le\infty$
be  the components of $F^*\cap I$. Denote by $\Ga_j$ the family of
all paths $\ga\subset\He$ which separate $F$ from $\infty$  such
that $\ov{\ga}\cap (c_j,d_j)\neq\emptyset$.

  Let $j'$ be such that
 $\tilde{v}_{j'}=\max_{3\le j\le N}\{\tilde{v}_j\}$ and let
  $r_{j'}:=e^{\tilde{v}_{j'}}.$
  By (\ref{2.3}), (\ref{3.nn4}) and the monotonicity of the
  capacity,
  we have
 \begin{eqnarray*}
 m(\Ga_{j'})&=& m(\Psi(\Ga_{j'}))\\
 &\ge&
 m(\{\ga_r=\{w\in\He:\, |w|=r\}:\, 1<r<r_{j'}\})\\
 &=&\frac{1}{\pi}\log r_{j'}
 \ge
 \frac{1}{\pi}\log\kap(K_F)=\frac{1}{\pi}\log\frac{\kap
 (I)}{\kap(F)},
 \end{eqnarray*}
 i.e.,
 \beq\label{4.1}
 \frac{\kap(F)}{\kap(I)}\ge\exp\{-\pi m(\Ga_{j'})\}.
 \eeq
 \begin{lem}\label{lem4.1}
 Under the above assumptions and definitions  the inequalities
 \beq\label{4.su1}
 e^{-\tilde{v}_{j'}}\le
 \frac{{\em\kap}(F)}{{\em\kap(I)}}\le\frac{4}{2+e^{\tilde{v}_{j'}}+e^{-\tilde{v}_{j'}}}
 \eeq
 and
  \beq\label{4.su2}
  \tilde{v}_{j'}\ge\pi(m(\Ga_{j})-1)
  \eeq
  hold for any $3\le j\le N$.

   Furthermore,
  if $\tilde{v}_{j'}\le\frac{\pi}{4}$, then for any $3\le j\le N$,
  \beq\label{4.4}
  \tilde{v}_{j'}\ge\frac{\pi}{2}\exp\left(-\frac{10\pi}{m(\Ga_{j})}\right).
  \eeq
  \end{lem}
  {\it Proof.}  The
  monotonicity of the capacity yields
 \begin{eqnarray*}
 r_{j'}&=&\kap(\{w:\,|w|\le r_{j'}\})\ge\kap(K_F)\\
 &\ge&\kap(\{\ov{\D}\cup[1,r_{j'}]\})
 =\frac{1}{4}\left(2+r_{j'}+\frac{1}{r_{j'}}\right),
 \end{eqnarray*}
 which, together with  (\ref{2.3}), implies (\ref{4.su1}).

 Next, we estimate from above the module of the family
 $\Ga'_{j}:=\psi(\Ga_{j})$.
 Consider
the metric
 $$\rho(w):=\left\{\begin{array}{ll}
 1,&\mb{ if }w\in\Si_F,\, \Im w\le \tilde{v}_{j'}+\pi,\\[2ex]
  0,&\mb{ elsewhere},
\end{array}\right.
 $$
 Since for any $\ga\in\Ga'_{j}$,
 $$\int_\ga\rho(w)|dw|\ge\pi,$$
 we see that
 $$
 m(\Ga'_{j})\le \frac{A(\rho)}{\pi^2}=
 \frac{\tilde{v}_{j'}+\pi}{\pi},
 $$
 from which (\ref{4.su2}) follows.

  For the
 small values of $\tilde{v}_{j'}$ we derive another estimate, i.e., inequality
 (\ref{4.4}), which
 reflects
the fact that $ m(\Ga_j)=m(\Ga'_{j})$ can be arbitrary small.
 Let $\tilde{v}_{j'}\le\frac{\pi}{4}.$ Without loss of generality we can assume that
 $\tilde{u}_{j}\ge\frac{\pi}{2}$ (if $\pi$ is closer to
 $\tilde{u}_{j}$ than $0$ the reasoning below has to be  modified
 in a straightforward way).

 Consider
the metric
 $$\rho_j(w):=\left\{\begin{array}{ll}
|w-\tilde{u}_{j}|^{-1},&\mb{ if }w\in\Si_F,\, \tilde{v}_{j'}\le
|w-\tilde{u}_{j}|\le
\pi,\\[2ex] 0,&\mb{ elsewhere},
\end{array}\right.$$
  We
claim that for any $\ga\in\Ga'_{j},$
 \beq\label{4.3}
\int_\ga\rho_j(w)\,
|dw|\ge\frac{1}{\sqrt{5}}\log\frac{\pi}{2\tilde{v}_{j'}}\,
 .\eeq
 In order to prove (\ref{4.3}), for $\tilde{v}_{j'}\le
r<R\le \tilde{u}_{j}$ we set
 \begin{eqnarray*}
B_j(r,R)&:=&\Si_F\bigcap\Biggl(\left\{w=u+iv:\, r\le
|u-\tilde{u}_{j}|\le
R,\, 0\le v\le \tilde{v}_{j'}\right\}\\
 &\bigcup& \Bigl\{ w=\tilde{u}_j+t
e^{i\theta}:\, r^2+\tilde{v}_{j'}^2\le t^2\le R^2+\tilde{v}_{j'}^2,\, \\
&&\hskip2cm
\sin^{-1}\left(\frac{\tilde{v}_{j'}}{t}\right)\le\theta\le\pi
-\sin^{-1}\left(\frac{\tilde{v}_{j'}}{t}\right) \Bigr\} \Biggr).
\end{eqnarray*}
For $\ga\in\Ga'_{j}$ and $R\le 2r$ we have
 $$\int_{\ga\cap B_j(r,R)}\rho_j(w)\,
|dw|\ge\frac{R-r}{\sqrt{\tilde{v}_{j'}^2+R^2}}\ge\frac{R-r}{\sqrt{5}r}\ge
\frac{1}{\sqrt{5}}\log\frac{R}{r}\, ,$$
 from which (\ref{4.3})
immediately follows.

According to the definition of the module (\ref{aaa}) and our
assumption that $\tilde{v}_{j'}\le\frac{\pi}{4}$ we obtain
 \begin{eqnarray*}
m(\Ga_{j}')&\le& 5\left(\log\frac{\pi}{2\tilde{v}_{j'}}\right)^{-2} A(\rho_j)\\
&\le& 5\,\pi
\left(\log\frac{\pi}{2\tilde{v}_{j'}}\right)^{-2}\log\frac{\pi}{\tilde{v}_{j'}}\le
 \frac{10\pi}{\log \frac{\pi}{2\tilde{v}_{j'}}},
\end{eqnarray*}
 which proves (\ref{4.4}).

 \hfill\rbox

 We use Lemma \ref{lemsu} below in Section 6.
 \begin{lem}\label{lemsu}
 Let $a<b$ and let $F\subset [a,b]$ be a regular compact set such
 that $a,b\in F$ and
 \beq\label{4.pi1}
 \frac{{\em\kap}(F)}{{\em\kap}([a,b])}\le e^{-3\pi}.
 \eeq
 Then, there exists a compact set $F^\star$ consisting of a finite
 number of closed intervals, such that $F\subset F^\star\subset [a,b]$ and
 \beq\label{4.pi2}
 0<q'<\frac{{\em\kap}(F^\star)}{{\em\kap}([a,b])}<q''<1
 \eeq
 holds with absolute constants $q'$ and $q''$.
 \end{lem}
 {\it Proof.} We can certainly assume that
 $[a,b]=I$.
 According to (\ref{4.pi1}) and the left-hand side of (\ref{4.su1})
  we have $\tilde{v}_{j'}\ge 3\pi$. Consider the set
 $$F^\star:=\{x\in I:\,\Im\psi(x)\le\tilde{v}_{j'}-2\pi\}.
 $$
 Let $\psi^\star:\He\to\Si_{F^\star},\{\tilde{u}_k^\star\}$ and
 $\{\tilde{v}_k^\star\}$ be defined as in Section 2 for the
 compact set $F^\star$ instead of $F$ and let
 $$\tilde{v}^\star_{k'}:=\max_k\{\tilde{v}_k^\star\}.$$
 Denote by $\Ga_k^\star$ the family of all paths $\ga\subset\He$
 which separate $F^\star$ from $\infty$ such that
 $\ov{\ga}\cap\psi^{\star-1}((\tilde{u}_k^\star,\tilde{u}_k^\star+
 i\tilde{v}_k^\star])\neq\emptyset$.
 Let $k^\star$ be defined such that
 $$
 [\tilde{u}_{k^\star}^\star,\tilde{u}_{k^\star}^\star+
 i\tilde{v}_{k^\star}^\star]:=\psi^\star\circ\psi^{-1}(
 [\tilde{u}_{j'}+i(\tilde{v}_{j'}-2\pi),\tilde{u}_{j'}+i\tilde{v}_{j'}]).
 $$
 Since
 \begin{eqnarray*}
 m(\Ga^\star_{k^\star})&=&m(\psi(\Ga^\star_{k^\star}))\\
 &\ge& m(\{[it,\pi+it]:\,
 \tilde{v}_{j'}-2\pi<t<\tilde{v}_{j'}\})=2,
 \end{eqnarray*}
 (\ref{4.su2}) written for $F^\star$ instead of $F$ shows that
 $$\tilde{v}_{k'}^\star\ge\pi.
 $$
 Therefore, the right-hand side of (\ref{4.su1}) implies the
 right-hand side of (\ref{4.pi2}).

 In order to prove the left-hand side of (\ref{4.pi2}), we consider
 the metric
 $$ \rho(w)=
\left\{\begin{array}{ll} 1\, ,& \mb{ if }0\le\Re
w\le\pi,\tilde{v}_{j'}-3\pi\le\Im w\le \tilde{v}_{j'}+\pi
,\\[2ex] 0,&\mb{
 elsewhere}.\end{array}\right.$$
  Since
\begin{eqnarray*}
\frac{\tilde{v}^\star_{k'}}{\pi}&=&m(\{[it,\pi+it]:\,
0<t<\tilde{v}^\star_{k'}\})\le m(\psi^\star(\Ga^\star_{k'}))\\
&=& m(\Ga^\star_{k'})=m(\psi(\Ga^\star_{k'}))
\le\frac{A(\rho)}{\pi^2}=4,
\end{eqnarray*}
the left-hand side of (\ref{4.su1}) implies the left-hand side of
(\ref{4.pi2}).

\hfill\rbox

  Next, let $E\subset\R,\phi(\cdot)=\phi(\cdot,E), E^*=\R\setminus
  E=\cup_j(c_j,d_j)$, and $[u_j,u_j+iv_j]=\phi([c_j,d_j])$ be defined as
  in Section 2. Consider $\al,\be\in\R\setminus(\cup_ju_j)$ such that
  $\al<\be$. The points
  $$a:=\phi^{-1}(\al),\quad b:=\phi^{-1}(\be)$$
  are uniquely defined. Let
  $$v_{[\al,\be]}:=\sup\{v_j:\, \al<u_j<\be\}.$$
 \begin{lem}\label{lem4.2}
 (a) If
 \beq\label{4.5}
 \frac{v_{[\al,\be]}}{\be-\al}\le q_1,
 \eeq
 then there exists $0<q_2=q_2(q_1)<1$ such that
 \beq\label{4.6}
 \frac{{\em \kap}(E\cap[a,b])}{{\em\kap}([a,b])}\ge q_2.
 \eeq
 (b) If
 \beq\label{4.7}
 \frac{v_{[\al,\be]}}{\be-\al}\ge q_3>0,
 \eeq
 then there exists $0<q_4=q_4(q_3)<1$ such that
 \beq\label{4.8}
 \frac{{\em\kap}(E\cap[a,b])}{{\em\kap}([a,b])}\le q_4.
 \eeq
 (c) If
 \beq\label{4.y1}
 \frac{{\em\kap}(E\cap[a,b])}{{\em\kap}([a,b])}\le q_5<1
 \eeq
 then there exists $q_6=q_6(q_5)>0$ such that
 \beq\label{4.y2}
 \frac{v_{[\al,\be]}}{\be-\al}\ge q_6.
 \eeq
 (d) If
 \beq\label{4.yy1}
 \frac{{\em\kap}(E\cap[a,b])}{{\em\kap}([a,b])}\ge q_7>0
 \eeq
 then there exists $q_8=q_8(q_7)>0$ such that
 \beq\label{4.yy2}
 \frac{v_{[\al,\be]}}{\be-\al}\le q_8.
 \eeq
 \end{lem}
 {\it Proof.} Let $F=h(E\cap[a,b])$, where
 \beq\label{4.1n}
 h(z)=\frac{2z-(a+b)}{b-a},
 \eeq
 and let $F^*\cap I=\cup_j(c_j,d_j),\psi:\He\to\Si_F$ and
 $\{\tilde{u}_j\}$, $\{\tilde{v}_j\}$ be defined as in Section 2
 with changed, for the convenience, index numbering to be the same as the
 numbering
 of $u_j$'s and $v_j$'s, i.e.,
 $$\psi\circ h\circ\phi^{-1}([u_j,u_j+iv_j])
 =[\tilde{u}_j,\tilde{u}_j+i\tilde{v}_j].$$
 Let $j'$ be such that $\tilde{v}_{j'}=\max_j\{\tilde{v}_{j}\}$
 and let $j^*$ be such that
 $v_{j^*}=v_{[\al,\be]},\al<u_{j^*}<\be$. As before, we denote by
 $\Ga_j$ the family of all paths $\ga\subset\He$ which separate
 $F$ from $\infty$ such that $\ov{\ga}\cap(c_j,d_j)\neq\emptyset$.
 Let $\tilde{\Ga}_j$ be the family of all half-ellipses in $\He$
 with the foci at $\al$ and $\be$ which have nonempty intersection
 with $(u_j,u_j+iv_j]$. By (\ref{3.nn3}) we have
 \beq\label{4.pi17}
 m(\Ga_j)=m(\phi\circ h^{-1}(\Ga_j))\ge
 m(\tilde{\Ga}_j)\ge\frac{1}{2\pi}\log\left(1+\frac{4v_j^2}{(\be-\al)^2}\right).
 \eeq
 Let $\Ga_j^*$ consist of all paths $\ga\subset\{w:\, \al<\Re
 w<\be,\Im w>0\}$ which separate $(\al,\be)$ from $\infty$
 such that $\ov{\ga}\cap(u_j,u_j+iv_j]\neq\emptyset.$ By the
 comparison principle we have
 \beq\label{4.pi7}
 m(\Ga_j)= m(\phi\circ h^{-1}(\Ga_j))\le m(\Ga_j^*).
 \eeq
 Note also that
 \beq\label{4.pi11}
 m(\Ga_{j'})=m(\psi(\Ga_{j'}))\ge m(\{(it,\pi+it):\,0<t<\tilde{v}_{j'}\})=
 \frac{\tilde{v}_{j'}}{\pi}.
 \eeq
 Next, applying the linear transformation $w\to\frac{\pi(w-\al)}{\be-\al}$
  we state the  analogues of (\ref{4.pi11}), (\ref{4.su2}), and
 (\ref{4.4}) for the module of the paths family $\Ga_j^*$:
 \beq\label{4.pi8}
 m(\Ga^*_{j^*})\ge m(\{(\al+it,\be+it):\,0<t<v_{[\al,\be]}\})
 =\frac{v_{[\al,\be]}}{\be-\al},
 \eeq
 \beq\label{4.pi9}
 m(\Ga^*_{j})\le\frac{v_{[\al,\be]}}{\be-\al}+1,
 \eeq
 and if $\frac{v_{[\al,\be]}}{\be-\al}\le\frac{1}{4}$ then
 \beq\label{4.pi10}
 m(\Ga_j^*)\le\frac{10\pi}{\log\frac{\be-\al}{2v_{[\al,\be]}}}.
 \eeq

 (a) By (\ref{4.5}), (\ref{4.pi7}) and (\ref{4.pi9})
 $$m(\Ga_{j'})\le q_1+1.$$
 By (\ref{4.pi11})
 $$
 \tilde{v}_{j'}\le\pi(q_1+1).
 $$
 Thus, applying the left-hand side of (\ref{4.su1}) and using the linear transformation
 (\ref{4.1n}) we obtain
 (\ref{4.6}).

 (b) By (\ref{4.7}) and (\ref{4.pi17})
 $$m(\Ga_{j^*})\ge C_1=C_1(q_3).
 $$
 Since by (\ref{4.4})
 $$\tilde{v}_{j'}\ge C_2=C_2(C_1),
 $$
   the right-hand side of (\ref{4.su1}) implies
  (\ref{4.8}).

  (c) By (\ref{4.y1}) and the left-hand side of (\ref{4.su1})
  $$ \tilde{v}_{j'}\ge \log\frac{1}{q_5}.
 $$
 Since by (\ref{4.pi7}) and (\ref{4.pi11})
 $$
 m(\Ga^*_{j'})\ge\frac{1}{\pi}\log\frac{1}{q_5},
 $$
 (\ref{4.pi10}) yields (\ref{4.y2}).

 (d) By (\ref{4.yy1}) and the right-hand side of (\ref{4.su1})
  $$ \tilde{v}_{j'}\le C_3=C_3(q_7).
 $$
 Since by (\ref{4.su2})
 $$
 m(\Ga_{j^*})\le\frac{C_3}{\pi}+1,
 $$
  (\ref{4.pi17})  implies (\ref{4.yy2}).

 \hfill\rbox

 \begin{lem}\label{lem4.3}
 Let $\al\be>0$ and, consequently, $ab>0$. Then
 \beq\label{4.12}
 \left(\frac{b-a}{b+a}\right)^2\le
 C\frac{(\be-\al)^2+v_{[\al,\be]}^2}{(\be+\al)^2}
 \eeq
 holds with the constant $C=2^{10}$.
 \end{lem}
 {\it Proof.} Without loss of generality we assume that $\be>\al>0$.
 Let $\Ga=\Ga(a,b,\He)$ be the family of all crosscuts of $\He$
 which separate points $a$ and $b$ from points $0$ and $\infty$. By
 (\ref{3.nn1})
 \beq\label{4.13}
 m(\Ga)\le\frac{1}{\pi}\log \frac{16 b}{b-a}.
 \eeq
 Let $\Ga_1$ be the family of all crosscuts of $\He_E$ which
 separate (in $\He_E$) the rectangle
 $$\{w:\, \al\le\Re w\le\be,0\le\Im w\le v_{[\al,\be]}\}
 $$
 from $0$ and $\infty$. By virtue of (\ref{3.nn2}) we have
 \beq\label{4.14}
 m(\Ga_1)\ge m(\Ga_9(\al,\be,v_{[\al,\be]}))\ge\frac{1}{2\pi}
 \log \frac{\be^2}{(\be-\al)^2+v_{[\al,\be]}^2}.
 \eeq
 Since
 $$m(\Ga)=m(\phi(\Ga))\ge m(\Ga_1),
 $$
 comparing (\ref{4.13}) and (\ref{4.14}) with the last estimate we
 obtain (\ref{4.12}):
 \begin{eqnarray*}
 \left(\frac{b-a}{b+a}\right)^2&<& \left(\frac{b-a}{b}\right)^2\le
 16^2 \frac{(\be-\al)^2+v_{[\al,\be]}^2}{\be^2}\\
 &\le& 2^{10} \frac{(\be-\al)^2+v_{[\al,\be]}^2}{(\be+\al)^2}.
 \end{eqnarray*}

 \hfill\rbox

 \begin{lem}\label{lem4.4}
 Let $E$ be such that
 \beq\label{4.15}
 \lim_{|u_j|\to\infty}\frac{v_j}{u_j}=0.
 \eeq
 Then there exist points $e_k,f_k,g_l,h_l\in E$ such that
 $\phi(e_k),\phi(f_k),\phi(g_k),\phi(h_k)\in\R\setminus\cup_ju_j$,
 $$\ldots\le e_{-1}<f_{-1}\le 0\le e_0<f_0\le e_1<f_1\le\ldots,$$
  $$\ldots\le g_{-1}<h_{-1}\le 0\le g_0<h_0\le g_1<h_1\le\ldots,$$
  $$E^*\subset
  \left[\cup_k(e_k,f_k)\right]\cup\left[\cup_l(g_l,h_l)\right],\quad
 (e_k,f_k)\cap (g_l,h_l)=\emptyset,$$
 $$
 \inf_k\frac{{\em\kap}(E\cap[e_k,f_k])}{{\em\kap}([e_k,f_k])}\ge q_1,
 $$
 \beq\label{4.17}
 \sum_k\left(\frac{f_k-e_k}{|e_k|+1}\right)^2<\infty,
 \eeq
 \begin{eqnarray}
 \frac{1}{6}&\le&\inf_l\frac{\max\{v_j:\,
 \phi(g_l)<u_j<\phi(h_l)\}}{\phi(h_l)-\phi(g_l)}\nonumber
 \\\label{4.hu18}
 & \le& \sup_l\frac{\max\{v_j:\,
 \phi(g_l)<u_j<\phi(h_l)\}}{\phi(h_l)-\phi(g_l)}\le
 \frac{1}{2},
 \end{eqnarray}
 \beq\label{4.18}
 q_2\le\inf_l\frac{{\em\kap}(E\cap[g_l,h_l])}{{\em\kap}([g_l,h_l])}
 \le \sup_l\frac{{\em\kap}(E\cap[g_l,h_l])}{{\em\kap}([g_l,h_l])}\le q_3,
 \eeq
 where $q_1>0$ and $0<q_2<q_3<1$ are absolute constants.
 \end{lem}
 {\it Proof.} Recall that by our assumption
 $0=\phi(0)$ belongs to $E$ with a certain closed interval around the origin.
  First, we introduce points
 $\al_m\in\R\setminus\cup_j u_j $ such that
 $$\ldots<\al_{-2}<\al_{-1}<\al_0=0<\al_1<\al_2\ldots,\quad
 1\le\al_{m+1}-\al_m<2.
 $$
 Let
 $$\de_{k,m}:=\max\{v_j:\, \al_k<u_j<\al_{m}\}\quad (k<m).
 $$
 If $\phi^{-1}([\al_k,\al_m])\subset E$ we set $\de_{k,m}:=0$.
 Next, we define a subsequence $\{\ga_s\}$ of a sequence
 $\{\al_m\}$ in the following way.

 Let $\ga_0:=\al_0=0$. If $\de_{0,1}\le\frac{1}{2}(\al_1-\al_0)$ we
 set $\ga_1:=\al_1$. Otherwise, i.e., if
 $\de_{0,1}>\frac{1}{2}(\al_1-\al_0)$,
 we consider the interval $(\al_0,\al_2)$. If
 $$\de_{0,2}\le\frac{1}{2}(\al_2-\al_0 )$$
 we set $\ga_1:=\al_2$. Otherwise, i.e., if
 $$\de_{0,2}>\frac{1}{2}(\al_2-\al_0 ),$$
 we consider the interval $(\al_0,\al_3)$, etc. According to our
 assumption (\ref{4.15}) after a finite number of steps we
 have
 $$\de_{0,m_1-1}>\frac{1}{2}(\al_{m_1-1}-\al_0 )
  \ge\frac{1}{6}(\al_{m_1}-\al_0 ),$$
  $$\de_{0,m_1}
  \le\frac{1}{2}(\al_{m_1}-\al_0 ).$$
 We set $\ga_1:=\al_{m_1}$ and proceed in the same way to construct
 $\ga_2:=\al_{m_2},m_2>m_1$ such that either $m_2=m_1+1$ and
 $$\de_{m_1,m_2}\le \frac{1}{2}(\al_{m_2}-\al_{m_1})$$
 or $m_2>m_1+1$ and
 $$\de_{m_1,m_2-1}>\frac{1}{2}(\al_{m_2-1}-\al_{m_1})
  \ge\frac{1}{6}(\al_{m_2}-\al_{m_1}),$$
  $$\de_{m_1,m_2}
  \le\frac{1}{2}(\al_{m_2}-\al_{m_1}).$$
  Repeating this procedure for both positive and negative indices
  $m$ we obtain the sequence of real numbers
  $$\ldots<\ga_{-2}<\ga_{-1}<\ga_0=0<\ga_1<\ga_{2}<\ldots$$
  with the following properties.

  (a) Either $1\le\ga_{s+1}-\ga_s<2$ and
  \beq\label{4.19}
  \max\{v_j:\,\ga_{s}<u_j<\ga_{s+1}\}\le\frac{1}{2}(\ga_{s+1}-\ga_s);
  \eeq
  or

  (b)
  \beq\label{4.20}
  \frac{1}{6}(\ga_{s+1}-\ga_s)\le
  \max\{v_j:\,\ga_{s}<u_j<\ga_{s+1}\}\le\frac{1}{2}(\ga_{s+1}-\ga_s).
  \eeq
  We denote
   by $\mu_k$ and $\mu_{k+1}$ in the case (a) and by $\nu_l$ and $\nu_{l+1}$ in the case
   (b) the endpoints of the intervals $(\ga_s,\ga_{s+1})$ . Let
   $$e_k:=\phi^{-1}(\mu_k)\quad f_k:=\phi^{-1}(\mu_{k+1}),$$
   $$g_l:=\phi^{-1}(\nu_l)\quad h_l:=\phi^{-1}(\nu_{l+1}).$$
   By (\ref{4.5})-(\ref{4.6}) and (\ref{4.19}) for any $k$ under
   consideration,
   we have
   $$
   \frac{\kap(E\cap[e_k,f_k])}{\kap([e_k,f_k])}\ge q_1>0.
   $$

 Since by Lemma \ref{lem4.3} and (\ref{4.19})
 $$\left(\frac{f_k-e_k}{|e_k|+1}\right)^2\le\frac{C}{k^2}
 $$
 holds with some constant $C>0$,
 we obtain (\ref{4.17}).

 Furthermore, (\ref{4.20}) yields (\ref{4.hu18}) which, together
 with  parts (a) and (b) of Lemma \ref{lem4.2}, implies (\ref{4.18}).

 \hfill\rbox

 We state another result that can be proved in a similar
 manner.
 \begin{lem}\label{lem4.41}
 Let $F\subset I,\pm1\in F$ be such that
 \beq\label{4.f15}
 \lim_{\tilde{u}_j\to u_0}\frac{\tilde{v}_j}{\tilde{u}_j-u_0}=0,
 \eeq
 where $u_0:=\psi(0)$.
 Then there exist points $\tilde{e}_k,\tilde{f}_k,\tilde{g}_l$, and
  $\tilde{h}_l\in F$ such that
  $$-1\le\tilde{e}_{-1}<\tilde{f}_{-1}\le\tilde{e}_{-2}<
  \ldots<0<\ldots <\tilde{f}_1\le\tilde{e}_0<\tilde{f}_0\le
  1,$$
  $$-1\le\tilde{g}_{-1}<\tilde{h}_{-1}\le\tilde{g}_{-2}<
  \ldots<0<\ldots <\tilde{h}_1\le\tilde{g}_0<\tilde{h}_0\le
  1,$$
  $$F^*\cap I\subset
  \left[\cup_k(\tilde{e}_k,\tilde{f}_k)\right]\cup\left[\cup_l(\tilde{g}_l,\tilde{h}_l)\right],\quad
 (\tilde{e}_k,\tilde{f}_k)\cap (\tilde{g}_l,\tilde{h}_l)=\emptyset,$$
 $$
 \inf_k
 \frac{{\em\kap}(F\cap[\tilde{e}_k,\tilde{f}_k])}{{\em\kap}([\tilde{e}_k,\tilde{f}_k])}\ge q_1,
 $$
 $$
 \sum_k\left(\frac{\tilde{f}_k-\tilde{e}_k}{\tilde{e}_k}\right)^2<\infty,
 $$
 \begin{eqnarray*}
 q_2&\le&\inf_l\frac{{\em\kap}(E\cap[\tilde{g}_l,\tilde{h}_l])}{
 {\em\kap}([\tilde{g}_l,\tilde{h}_l])}\\
 &\le&\sup_l\frac{{\em\kap}(E\cap[\tilde{g}_l,\tilde{h}_l])}{
 {\em\kap}([\tilde{g}_l,\tilde{h}_l])}\le q_3,
 \end{eqnarray*}
 where $q_1>0$ and $0<q_2<q_3<1$ are absolute constants.
 \end{lem}

\absatz{Proof of Theorem \ref{th1}}

{\it Proof of (i)$\Rightarrow$(ii).}
 Let $E_R$ and $\phi_R$ be defined as in Section 2. Denote by $u_{j,R}$ and $v_{j,R}$
 the appropriate real numbers $u_j$ and $v_j$ defined in Section 2 for $E_R$ instead of $E$.
 By Theorem C we can assume that $[-e^{3\pi^2},e^{3\pi^2}]\subset E$
 and $a_j\ge e^{2\pi^2}$ for $j\ge 0$ as well as $b_j\le -e^{2\pi^2}$ for
 $j<0$.

 Since by (\ref{1.5})
 $$
 \lim_{j\to\pm\infty}\frac{b_j-a_j}{a_j}=0,
 $$
 we can also  assume that
 $$\frac{b_j-a_j}{\min\{|a_j|,|b_j|\}}<\frac{\pi}{2}.
 $$
 Let
 $$K_t:=\{z\in\He:\,|z|=t\}\quad (t>0),$$
 and let
 $$u^+_R:=\inf\{u_{j,R}:\, u_{j,R}>0\},$$
  $$u^-_R:=\sup\{u_{j,R}:\, u_{j,R}<0\}.$$
  We set $u^+_R:=\infty$ if $E^*\cap[0,R]=\emptyset$ and
  $u^-_R:=-\infty$ if $E^*\cap[-R,0]=\emptyset$. Let
  $$u_R:=\min\{u^+_R,-u^-_R\}<\infty,$$
  $$s_R:=\sup\{|\phi_R(w)|:\, w\in K_1\}.$$
  We start with the observation that the
  inequality
  \beq\label{5.p1}
  s_R\le u_R
  \eeq
  holds for sufficiently large $R>e^{3\pi^2}$.

  Indeed, let   $z^*_R\in \ov{K_1}$ be such that for
  $w^*_R:=\phi_R(z^*_R)$ we have $|w^*_R|=s_R$.
  Assume,  contrary to our claim (\ref{5.p1}), that
  $$
  s_R=|w^*_R|>u_R.
  $$
  Let $\Ga_1=\{ K_t:\, 1<t<\exp(3\pi^2)\}$. By (\ref{3.nn4}) we
  have
  $$
  m(\Ga_1)=3\pi.
  $$
  Consider the family $\Ga_1'=\Ga_1'(R):=\phi_R(\Ga_1)$ and the
  metric
  $$\rho(w):=\left\{\begin{array}{ll}
 1,&\mb{ if }w\in\He, |w|\le 2u_R,\\[2ex]
 0,&\mb{ elsewhere}.
\end{array}\right.$$
 Since
 $$
  m(\Ga'_1)\le\frac{A(\rho)}{u_R^2}=2\pi,
  $$
  we have the contradiction
  $$
  3\pi=m(\Ga_1)=m(\Ga_1')\le 2\pi.
  $$
  Hence, our assumption is incorrect, which proves
 (\ref{5.p1}).

 In the reasoning below $R$ and $T>R$ are sufficiently large.
 By (\ref{2.2})  we have
 $$
 |\phi_R(z)|\le 2T\quad (z\in K_T).
$$
 We proceed to show that
 \beq\label{5.2}
 |\phi_R(z)|\ge C|w_R|\quad (z\in K_1),
 \eeq
 where $w_R:=\phi_R(i)$ and $C=\exp(-2\pi^2)$.

 To prove this, let $w:=\phi_R(z), z\in K_1$.
 The only nontrivial  case is where $|w|<|w_R|$.
  Let $\Ga'_2=\Ga'_2(w,w_R)$  be the family of all
 crosscuts of $\He_{E_R}$ which separate $0$ and $w$ from $w_R$
 and $\infty$. Then (\ref{3.nn4}) and (\ref{5.p1}) imply
 \beq\label{5.3}
 m(\Ga'_2)\ge m(\{K_t:\,
 |w|<t<|w_R|\})=\frac{1}{\pi}\log\frac{|w_R|}{|w|}.
 \eeq
  Considering the metric
 $$\rho(z):=\left\{\begin{array}{ll}
 1,&\mb{ if }z\in\He,\, |z|\le
 2,\\[2ex] 0,&\mb{ elsewhere},
\end{array}\right.$$
 for $\Ga_2=\phi^{-1}_R(\Ga'_2)$
  we obtain
 \beq\label{5.4}
 m(\Ga_2)\le A(\rho)=2\pi.
 \eeq
 Comparing (\ref{5.3}) and (\ref{5.4}) we have (\ref{5.2}).

 Let
 $$t_R:=\frac{1}{2}C|w_R|,\quad C=\exp(-2\pi^2).
 $$
 According to (\ref{3.nn5}), for the family of radial intervals
 $$\Ga'_3=\Ga'_3(R,T)=\{\ga_\theta:=\{it_R+re^{i\theta}:\,
 t_R<r<2T\}:\, 0<\theta<\pi\}
 $$
 we have
  \beq\label{5.5}
  m(\Ga'_3)\ge\pi\left(\log\frac{2T}{t_R}\right)^{-1}
  =\pi\left(\log\frac{4T}{C|w_R|}\right)^{-1}.
  \eeq
 Our next objective is to estimate the module of the path family
 $\Ga_3=\phi^{-1}_R(\Ga_3')$ from above. Consider the set
 $$\tilde{E}_R:=\{x\in\R:\, \Im \phi_R(x)\le t_R\}\supset E_R$$
 which consists of a finite number of closed intervals.

  Let $\{(a_j,b_j)\}_{j=-J^-}^{J^+}$  be the intervals from the part (i) of Theorem
  \ref{th1} satisfying the condition
  \beq\label{5.7}
  \R\setminus
  \tilde{E}_R=:\tilde{E}^*_R\subset\bigcup_{j=-J^-}^{J^+}(a_j,b_j)\subset(-T,T).
  \eeq
  We assume that the system of intervals
  $\{(a_j,b_j)\}_{j=-J^-}^{J^+}$ is minimal in the sense that
  $J^{\pm}$ cannot be decreased with (\ref{5.7}) still valid.
 It can happen that there are no such intervals at all (in this
 case we write $J^-=0$ and $J^+=-1$) or there are only intervals
 with positive endpoints (that is $J^-=0$) or there are only
 intervals with negative endpoints (that is $J^+=-1$).

 Let $\Ga_{13}=\Ga_{13}(T,\tilde{E}_R\cap([1,T]\cup[-T,-1])
 ,\{(a_j,b_j)\}_{j=-J^-}^{J^+})$ be the path
 family from Section 3. Notice that
 \beq\label{5.nn1}
 m(\Ga_3)\le m(\Ga_{13})
 \eeq
 and by (\ref{1.5}) and Lemma \ref{lem3.1n}
 \beq\label{5.n1}
 m(\Ga_{13})\le\frac{\pi(\log T +C_1)}{(\log T)^2}.
 \eeq
  Here and in the sequel we adopt the convention that $C,C_1,C_2,\ldots$
 denote positive constants, possibly different in different cases.

 Referring to  (\ref{5.5}),  (\ref{5.nn1}), and (\ref{5.n1}) we
 find that
 $$
 \left(\log T\right)^2\le\left(
 \log\frac{4T}{C|w_R|}\right)(\log T+C_1),
 $$
 and  making $T\to\infty$ we see  that
 $$
 \log|w_R|\le \log\frac{4}{C}+C_1.
 $$
 This means that $|w_R|=|\phi_R(i)|$ is uniformly bounded and we can apply Lemma
 \ref{lem2.2} to derive (ii).

 \hfill\rbox

 {\it Proof of (ii)$\Rightarrow$(i).}
 We first  show that (\ref{4.15}) holds.
 Let $\Ga=\Ga(y),y>1$ be the family of all  curves in $\He$
 which separate $i$ and $iy$ from $\R$. By (\ref{3.1}) we have
 \beq\label{5.10}
 m(\Ga)\ge\frac{\pi}{4\log 4y}.
 \eeq
 Let $\Ga':=\phi(\Ga)$,
 $w_0:=\phi(i)$, and $w_1:=\phi(iy).$ We assume that $y$ is
 sufficiently large. In particular, $y$ is so large  that  $|w_0|<|w_1|$. Let
 $$Q=Q(|w_0|,|w_1|):=\{w\in\He:\, |w_0|\le|w|\le|w_1|\}.
 $$
 Denote by $\Ga_1'$ the family of all curves in $\He_E\cap
 Q$ joining circular arcs $\{w:|w|=|w_0|\}$ and
 $\{w:|w|=|w_1|\}$. Since any $\ga\in\Ga'$ includes two disjoint
 curves from the family $\Ga_1'$, we conclude that
 \beq\label{5.11}
 m(\Ga')\le\frac{m(\Ga_1')}{4}.
 \eeq
 We prove (\ref{4.15}) by contradiction. Suppose
  it were false. Then, we could find  a constant $0<c<1$ and a
 (monotone) sequence of integers $\{j_k\}_{k=-M}^N$,
 where $M+N=\infty$, such that
 $$v_{j_k}\ge c|u_{j_k}|,$$
 $$u_{j_{k+1}}>2u_{j_k}>2|w_0|\quad (j_k\ge0),$$
 $$u_{j_{k-1}}<2u_{j_k}<-2|w_0|\quad (j_k<0).$$
 Consider the domain
 $$D=D(\{u_{j_k}\},c, |w_1|):=Q\setminus
 \bigcup_{k=-M_1}^{N_1}[u_{j_k},u_{j_k}(1+ic)],$$
 where $M_1\le M$ and $N_1\le N$ are such that for $0\le k\le
 N_1$ or $-M_1\le k\le -1$ we have
 $$|u_{j_k}|<\frac{|w_1|}{2}.$$
 Notice that
  \beq\label{5.nn3}
  m(\Ga_1')\le m(\Ga_2'),
  \eeq
  where $\Ga_2'$ is the family of all curves in $D$ joining
  circular arcs $\{w:\, |w|=|w_0|\}$ and $\{w:\, |w|=|w_1|\}$.

 Applying  Lemma \ref{lema} with $r=|w_0|$ and $R=\frac{|w_1|}{2}$
 we obtain
 \beq\label{5.nn4}
 m(\Ga_2')\le\frac{\pi\log\frac{|w_1|}{|w_0|}-C(M_1+N_1)
 }{\left(\log\frac{|w_1|}{|w_0|}\right)^2}.
 \eeq
 According to Lemma  \ref{lem2.1} we can choose $y$ to be
 arbitrarily large and to satisfy the inequality
 \beq\label{5.sa1}
 y\le C_1|w_1|.
 \eeq
 Therefore, comparing (\ref{5.10})-(\ref{5.nn4})
  for such values of $y$ we have
 $$
 \frac{\pi}{\log 4C_1|w_1|}\le \frac{\pi\log C_2|w_1|-C(M_1+N_1)}
 {\left(\log C_2|w_1|\right)^2},
 $$
 i.e.,
 $$
 M_1+N_1\le \frac{C_3\log C_2|w_1|}{\log C_1|w_1|}.
 $$
 Passing to the limit as  $y\to\infty$, i.e., as $|w_1|\to\infty$ we
 have $M+N\le C_3$. This contradiction proves (\ref{4.15}).

 As a cover $\{(a_j,b_j)\}$ of $E^*$ satisfying (\ref{1.41n})-(\ref{1.5})
  we  use the cover $\{(e_k,f_k)\}$,
 $\{(g_l,h_l)\}$, constructed for any $E$ satisfying (\ref{4.15})
  in Lemma \ref{lem4.4}. Thus,
 we only need  to show that
 \beq\label{5.9}
 \sum_l\left(\frac{h_l-g_l}{|g_l|+1}\right)^2<\infty.
 \eeq
 Let $\al_l:=\phi(g_l),\be_l=\phi(h_l)$. According to
 (\ref{4.hu18}) we have
 \beq\label{5.13}
 \frac{1}{6}(\be_l-\al_l)\le v'_l:=\max\{v_j:\,
 \al_l<u_j<\be_l\}\le \frac{1}{2}(\be_l-\al_l).
 \eeq
 Let $u'_l\in (\al_l,\be_l)$ be such that
 $[u'_l,u_l'+iv_l']\subset\partial\He_E$ and
 let $l^\pm\le L^\pm$ be defined such that for $l^+\le l\le L^+$
 and $-L^-\le l\le -l^-$,
 $$(\al_l,\be_l)\subset \left[|w_0|,\frac{|w_1|}{2}
 \right]\cup \left[-\frac{|w_1|}{2},-|w_0|\right]
 $$
 and $v_l'<|u_l'|$.

 Consider the domain
 $$D_1=D_1(\{u'_{l}\},\{v'_{l}\}, |w_0|, |w_1|):=Q\setminus
 \bigcup_{l^+\le l\le L^+\atop -L^-\le l\le -l^-}[u'_{l},u'_{l}+iv'_l].$$
 Notice that
  \beq\label{5.nnn5}
  m(\Ga_1')\le m(\Ga_3'),
  \eeq
  where $\Ga_3'$ is the family of all curves in $D_1$ joining
  circular arcs $\{w:\, |w|=|w_0|\}$ and $\{w:\, |w|=|w_1|\}$.

 By Lemma \ref{lema}
 we have
 \beq \label{5.nnn6}
 m(\Ga_3')\le\frac{\pi\log\frac{|w_1|}{|w_0|}-C_4\sum_{l^+\le l\le L^+\atop -L^-\le l\le -l^-}
 \left(\frac{v'_l}{u'_l}\right)^2}{\left(\log\frac{|w_1|}{|w_0|}\right)^2}.
 \eeq
 Comparing (\ref{5.10}), (\ref{5.11}),(\ref{5.nnn5}), and
 (\ref{5.nnn6}) for $y$ satisfying (\ref{5.sa1})
 we have
 $$ \frac{1}{\log 4C_1|w_1|} \le
 \left(\log C_2|w_1|\right)^{-2}\left(
 \log C_2|w_1|-\frac{C_4}{\pi}\sum_{l^+\le l\le L^+\atop
  -L^-\le l\le -l^-}\left(\frac{v'_l}{u'_l}\right)^2\right),
  $$
  i.e.,
  $$
 \sum_{l^+\le l\le L^+\atop
  -L^-\le l\le
  -l^-}\left(\frac{v'_l}{u'_l}\right)^2\le C_5\frac{\log
  C_2|w_1|}{\log 4C_1|w_1|}.
  $$
    Passing to the limit as $y\to\infty$ and applying (\ref{5.13}) we obtain
  $$ \sum_l\frac{(\be_l-\al_l)^2+v_l'^2}{(\be_l+\al_l)^2}\le C_6
  \sum_{l}\left(\frac{v'_l}{u'_l}\right)^2\le C_7.
  $$
  The last inequality and  Lemma \ref{lem4.3} imply
   (\ref{5.9}).

 \hfill\rbox

 As we mentioned in the introduction, we can reformulate the part
 (i)$\Leftrightarrow$(iii) of Theorem \ref{th1} in the form of the
 equivalence (i')$\Leftrightarrow$(iii') of Remark 5.

 {\it Proof of (i')$\Rightarrow$(iii') of Remark 5.} We use the same idea
 as in the proof of the part (i)$\Rightarrow$(ii) of
 Theorem \ref{th1}.

 Let $u_0:=\psi(0),$ and let for $r>0$,
 $$K_r:=\{z\in\He:\, |z|=r\},\quad K'_r:=\psi(K_r).$$
 Denote by $w_r\in K_r'$ any point with the property
 $\Re w_r=u_0$ and let $z_r:=\psi^{-1}(w_r).$

 Note that for any $w\in K'_r$,
 \beq\label{5.15}
 |w-u_0|\le C|w_r-u_0|,\quad C=\exp(2\pi^2).
 \eeq
 Indeed, if $|w-u_0|>|w_r-u_0|$, we consider the family $\Ga=\Ga(r)$ of
 all crosscuts of $\He$ which separate $0$ and $z_r$ from $z:=\psi^{-1}(w)$ and
 $\infty$.
 Considering the metric
 $$ \rho(\z)=
\left\{\begin{array}{ll} 1\, ,& \mb{ if }\z\in
 \He, |\z|\le 2r,\\[2ex] 0,&\mb{
otherwise},\end{array}\right.
 $$
 we have
 \beq\label{5.51}
 m(\Ga)\le \frac{A(\rho)}{r^2}=2\pi.
 \eeq
 On the other hand, by (\ref{3.nn4}) for the family $\Ga'=\psi(\Ga)$ we obtain
 \beq\label{5.52}
 m(\Ga')\ge m(\{K_t:\, |w_r-u_0|<t<|w-u_0|\})=\frac{1}{\pi}\log
 \frac{|w-u_0|}{|w_r-u_0|}.
 \eeq
 Comparing (\ref{5.51}) and (\ref{5.52}) we have (\ref{5.15}).

 According to (\ref{5.15}),
 in order to verify (\ref{1.15}),  it is enough to show that
 \beq\label{5.16}
 |w_r-u_0|\le C_1r
 \eeq
 holds with a constant $C_1>0$ independent of $r$
 for sufficiently small values of $r>0$. In particular, we can assume
 that $|w_r-u_0|<\frac{\log2}{2}$.

 Consider the family $\Ga_1'=\Ga_8(a,b,c)$ from Section 3 with
 $a=u_0,b=\frac{|w_r-u_0|}{3}$ and $c=\Im w_r=|w_r-u_0|$. By (\ref{3.51})
 \beq\label{5.53}
 m(\Ga_1')\ge\frac{\pi}{4\log\frac{18}{|w_r-u_0|}}.
 \eeq
 Let
 $$\tilde{F}_r:=\left\{x\in I:\, \Im
 \psi(x)\le\frac{|w_r-u_0|}{3}\right\}\supset F.
 $$
 Note that $\tilde{F}_r$ consists of a finite number of closed
 intervals.

 Denote by $\Ga_2=\Ga_2(r)$  the family of all paths in a half-ring
 $\{z\in\He:\, r<|z|<1\}$ which separate either $\tilde{F}_r\cap[r,1]$
 from $[-1,-r]$ or
 $\tilde{F}_r\cap[-1,-r]$ from $[r,1]$.
 The conformal invariance of the module and the comparison principle
 imply
 \beq\label{5.54}
 m(\Ga_1')\le \frac{1}{4} m(\psi(\Ga_2))=\frac{1}{4}m(\Ga_2).
 \eeq
 Let $\Ga_3$ consist of all $\ga\in\Ga_3$ possessing the following
 property: there exists $\ga^*\in\Ga_2$ such that $\ga=\{ z:\,
 \frac{1}{z}\in\ga^*\}$.

 Applying Lemma \ref{lem3.1n} and the assumptions (\ref{1.13})-(\ref{1.14})
 in the same way as we did in the proof of (\ref{5.n1}) we see that
 \beq\label{5.17}
 m(\Ga_2)=m(\Ga_3)\le\pi\left(\log\frac{C_2}{r}\right)^{-2}\log\frac{C_3}{r}.
 \eeq
 Comparing   (\ref{5.53})-(\ref{5.17})  we obtain
 (\ref{5.16}).

  \hfill\rbox

  {\it Proof of (iii')$\Rightarrow$(i') of Remark 5.} We use essentially the
  same idea as in the proof of the part (ii)$\Rightarrow$(i) of Theorem
  \ref{th1}. Let $\tilde{u}_j$ and $\tilde{v}_j$ be defined as in
 Section 2.
  By \cite[Corollary 1.12]{cartot} and
  \cite[Theorem  1]{and04}
   we have
 $$
 \lim_{t\to0^+}\frac{\kap(F\cap[0,t])}{\kap([0,t])}=
 \lim_{t\to0^+}\frac{\kap(F\cap[-t,0])}{\kap([-t,0])}=1,
 $$
 which implies (\ref{4.f15}).

 To see this, let for $\al\in(u_0,\pi)\setminus\cup_j\tilde{u}_j$,
 $$\de(\al):=\max\{\tilde{v}_j:\, u_0<\tilde{u}_j<\al\},$$
 and let $j^*$ be such that $u_0<\tilde{u}_{j^*}<\al$ and
 $v_{j^*}=\de(\al).$ Let $t:=\psi^{-1}(\al)$. Denote by
 $\Ga_\al$ the family of all paths $\ga\subset\He$ which
 separate $F\cap[0,t]$ from $\infty $  such that
 $\ov{\ga}\cap\psi^{-1}((u_{j^*},u_{j^*}+iv_{j^*}])\neq\emptyset$.
 By Lemma \ref{lem4.1}, applied to the set
 $\{x\in I:\,\frac{t}{2}(x+1)\in F\cap[0,t]\}$ instead of $F$, we obtain
 $$\frac{\kap(F\cap[0,t])}{\kap([0,t])}\le\frac{4}{2+e^{v(t)}+e^{-v(t)}},$$
 where
 $$v(t)\ge
\left\{\begin{array}{ll} \pi(m(\Ga_\al)-1),& \mb{ if
}v(t)>\frac{\pi}{4},
 \\[2ex] \frac{\pi}{2}\exp\left(-\frac{10\pi}{m(\Ga_\al)}\right),&\mb{
if }v(t)\le\frac{\pi}{4}.\end{array}\right.$$
 Passing to the limit as $\al\to u_0^+$, i.e., as $t\to 0^+$, we have
 $$\lim_{t\to 0^+} v(t)=0,
 $$
 that is,
 \beq\label{4.p1}
 \lim_{\al\to u_0^+} m(\Ga_\al)=0.
 \eeq
 Let $\Ga'_{1}=\Ga'_1(t)$ be the family of all half-ellipses in $\He$ with
 the foci at $u_0$ and $\al$ which have nonempty intersection with the
 vertical interval $(\tilde{u}_{j^*},\tilde{u}_{j^*}+i\tilde{v}_{j^*})$. Since
 \beq\label{5.hu1}
 m(\Ga_\al)=m(\psi(\Ga_\al))\ge m(\Ga'_{1}),
 \eeq
 according to (\ref{3.nn3}), (\ref{4.p1}) and (\ref{5.hu1}) we obtain
 \beq
 \label{fin.1}
 \lim_{\al\to u_0^+}\frac{\de(\al)}{\al-u_0}=0.
 \eeq
 The same conclusion can be drawn for
 $$\de(\al):=\max\{\tilde{v}_j:\, \al<\tilde{u}_j<u_0\}\quad
 0<\al<u_0,
 $$
 i.e.,
 \beq
 \label{fin.2}
 \lim_{\al\to u_0^-}\frac{\de(\al)}{u_0-\al}=0.
 \eeq
 The  relations (\ref{fin.1}) and (\ref{fin.2}) imply (\ref{4.f15}).

 As a cover $\{(a_j,b_j)\}$ of $F^*\cap I$ satisfying (\ref{1.13}) and
 (\ref{1.14}) we can use the cover $\{(\tilde{e}_k,\tilde{f}_k)\}$,
 $\{(\tilde{g}_l,\tilde{h}_l)\}$,  for any $F$ satisfying (\ref{4.f15})
 discussed
  in Lemma \ref{lem4.41}. In this case
 we only  need to show that
 \beq\label{5.91}
 \sum_l\left(\frac{\tilde{h}_l-\tilde{g}_l}{\tilde{g}_l}\right)^2<\infty.
 \eeq
 Let $\Ga_2=\Ga_2(y),0<y<1$ be the family of all  curves in $\He$
 which separate $i$ and $iy$ from $\R$. By (\ref{3.1}) we have
 \beq\label{5.s1}
 m(\Ga_2)\ge\frac{\pi}{4\log \frac{4}{y}}.
 \eeq
 Let $\Ga'_2:=\psi(\Ga_2)$,
 $w_0:=\psi(i)$, and $w_1:=\psi(iy).$ We assume that $y$ is
 sufficiently small. In particular, it is such that $2
 |w_1-u_0|<d_0:=\min\{u_0,\pi-u_0,|w_0-u_0|\}$. Let
 $$Q:=\{w\in\He:\, |w_1-u_0|\le|w-u_0|\le d_0\}
 $$
 and let $\tilde{\al}_l:=\psi(\tilde{g}_l),\tilde{\be}_l=\psi(\tilde{h}_l)$.
 Furthermore,
 let $\tilde{u}'_l\in (\tilde{\al}_l,\tilde{\be}_l)$ be such that
 $[\tilde{u}'_l,\tilde{u}_l'+i\tilde{v}_l']\subset\partial\Sigma_F$, and
 let positive integers $l^\pm\le L^\pm$ be defined so that for $l^+\le l\le L^+$
 and $-L^-\le l\le-l^-$,
 $$(\tilde{\al}_l,\tilde{\be}_l)\subset \left[u_0+|w_1-u_0|,u_0+\frac{d_0}{2}
 \right]\cup \left[u_0-\frac{d_0}{2},u_0-|w_1-u_0|\right]
 $$
 and $\tilde{v}_l\le|\tilde{u}_l-u_0|$.

 Consider the domain
 $$D=D(\{\tilde{u}'_{l}\},\{\tilde{v}'_{l}\}, |w_1-u_0|,d_0):=Q\setminus
 \bigcup_{l^+\le l\le L^+\atop -L^-\le l\le -l^-}[\tilde{u}'_{l},\tilde{u}'_{l}+i\tilde{v}'_l].$$
 Notice that
  \beq\label{5.s3}
  m(\Ga_2')\le \frac{1}{4}m(\Ga_3'),
  \eeq
  where $\Ga_3'$ is the family of all curves in $D$ joining
  circular arcs $\{w:\, |w-u_0|=|w_1-u_0|\}$ and $\{w:\, |w-u_0|=d_0\}$.

 By Lemma \ref{lema}
 we have
 \beq \label{5.s4}
 m(\Ga_3')\le\frac{\pi\log\frac{d_0}{|w_1-u_0|}-C\sum_{l^+\le l\le L^+\atop -L^-\le l\le -l^-}
 \left(\frac{\tilde{v}'_l}{\tilde{u}'_l-u_0}\right)^2}{\left(\log\frac{d_0}{|w_1-u_0|}\right)^2}.
 \eeq
 Further, proceeding  as in the proof of the part (ii)$\rightarrow$(i) of
Theorem \ref{th1} we derive  the inequality (\ref{5.91}) from
(\ref{5.s1})-(\ref{5.s4}). Since the proof of (\ref{5.91}) is
similar, we leave out the details of the proof.

 \hfill\rbox

 \absatz{Proof of Theorem \ref{th2}}

 The equivalence (ii)$\Leftrightarrow$(iii) is trivial (see
 Theorem C  and Theorem \ref{th1}).

 {\it Proof of (i)$\Rightarrow$(ii).}
 Let
 $$\al_j:=\phi(a_j),\, \be_j:=\phi(b_j),\, \delta_j:=\sup\{v_k:\,
 \al_j< u_k<\be_j\}.$$
 If necessary slightly moving  points $a_j$ and $b_j$,  we can assume that
 $\al_j,\be_j\in\R\setminus\cup_{j=-N}^Mu_j$.
 We also assume  that
 $\lim_{j\to \infty}\frac{\de_j}{\be_j}=0$ if $M=\infty$ and that $\lim_{j\to
 -\infty}\frac{\de_j}{\al_j}=0
 $ if $N=\infty$
 (because otherwise  (\ref{4.15}) does not hold and
 by the reasoning in the beginning  of the proof of the part
 (ii)$\Rightarrow$(i) of Theorem \ref{th1} we have dim $\cP_\infty\neq2,$ i.e.,
 by Theorem C dim $\cP_\infty=1$).

 Furthermore, extending $E\cap[a_j,b_j]$ for each interval
 $[a_j,b_j]$  with the property
 $$\frac{\kap(E\cap[a_j,b_j])}{\kap([a_j,b_j]}<e^{-3\pi}
 $$ in the way described in Lemma \ref{lemsu}, we obtain a new set
 $E^\star\supset E$ which satisfies (\ref{1.6}) in the strengthened
 form
 \beq\label{6.d1}
 0<q'\le\inf_j\frac{\kap(E^\star\cap[a_j,b_j])}{\kap([a_j,b_j]}
 \le\sup_j\frac{\kap(E^\star\cap[a_j,b_j])}{\kap([a_j,b_j]}\le
 q''<1,
 \eeq
 where $q'$ and $q''$ are constants.

 Since dim $\cP_\infty(\C\setminus E^\star)\ge$ dim $\cP_\infty(\C\setminus
 E)$, it is enough to show that \newline
  dim $\cP_\infty(\C\setminus
 E^\star)=1$. Hence, we can proceed  assuming that
 $E=E^\star$.

 The parts (c) and (d) of Lemma \ref{lem4.2} and (\ref{6.d1}) imply
 \beq\label{6.d2}
 C_1(\be_j-\al_j)\le\de_j\le C_2(\be_j-\al_j).
 \eeq
 Since by
 (\ref{1.7}) and (\ref{4.12}) we have
 $$\sum_{j=-N}^M\frac{(\be_j-\al_j)^2+\de_j^2}{(\be_j+\al_j)^2}=\infty,
 $$
  (\ref{6.d2})
 yields
 \beq\label{6.d3}
 \sum_{j=-N}^M\left(\frac{\de_j}{\al_j}\right)^2=\infty.
 \eeq
 We claim that
 \beq\label{6.si1}
 \lim_{y\to+\infty}\frac{|\phi(iy)|}{y}=0.
 \eeq
 Indeed, as in the proof of the part (ii)$\Rightarrow $(i) of
 Theorem \ref{th1} we let $w_0:=\phi(i)$ and $w_1:=\phi(iy),y>1$.
 We are interested in the case of sufficiently large values of
 $y$.
 Let $\Ga=\Ga(y)$ be the family of all  curves in $\He$
 which separate $i$ and $iy$ from $\R$. By (\ref{5.10}),
 (\ref{5.11}), (\ref{5.nnn5}), and (\ref{5.nnn6}) we have
 $$
 \frac{\pi}{4\log 4y}\le m(\Ga)\le
 \frac{\pi\log\frac{|w_1|}{|w_0|}-C_3\sum_{l^+\le l\le L^+\atop -L^-\le l\le -l^-}
 \left(\frac{\de_l}{\al_l}\right)^2}{4\left(\log\frac{|w_1|}{|w_0|}\right)^2},
 $$
 i.e.,
 $$
 \sum_{l^+\le l\le L^+\atop -L^-\le l\le -l^-}
 \left(\frac{\de_l}{\al_l}\right)^2\le\frac{\pi}{C_3}\frac{\log\frac{|w_1|}{|w_0|}}{\log
 4y}\log\frac{4|w_0|y}{|w_1|}.
 $$
 Comparing the last inequality with (\ref{6.d3}), we obtain
 (\ref{6.si1}) which, together with
 Lemma \ref{lem2.1}, implies   dim $\cP_\infty=1$.

  \hfill\rbox

 {\it Proof of (ii)$\Rightarrow$(i).}
 First, let (\ref{4.15}) hold.
 Consider the cover $\{(e_k,f_k)\}$ and $\{(g_l,h_l)\}$ of
 $E^*$ from Lemma \ref{lem4.4}.
 We can take the system of intervals $(g_l,h_l)$ to be our intervals
  $(a_j,b_j)$. By virtue  of (\ref{4.18})
 they satisfy (\ref{1.6}). They also satisfy (\ref{1.7}) because
 otherwise the whole system of intervals $\{(e_k,f_k)\}$ and
 $\{(g_l,h_l)\}$ satisfies (\ref{1.5}), which would mean by Theorem
 \ref{th1} that
 dim $\cP_\infty=2$, a contradiction.

 We now proceed with the case where
 $$\limsup_{|u_j|\to\infty}\frac{v_j}{|u_j|}>0.$$
 Let $u_{j_k}$,  $ v_{j_k}$ and $0<c<1$ be such that
 $$\lim_{k\to\infty}|u_{j_k}|=\infty,\quad v_{j_k}\ge c|u_{j_k}|,
 \quad v_{j_{k+1}}\ge v_{j_k}.$$
 There is no loss of generality in assuming that $u_{j_k}>0$ and $u_{j_{k+1}}>4u_{j_k}$
 (the case
 of the infinite number of negative $u_{j_k}$'s can be treated similarly).

 Let $\ga_{k}^\pm\in\R\setminus\cup_ju_j$ satisfy
 $$\ga_{k}^-< u_{j_k}<\ga_{k}^+,\quad \ga_{k}^+-\ga_{k}^-=u_{j_k},$$
and let
 $a_k:=\phi^{-1}(\ga_{k}^-),b_k:=\phi^{-1}(\ga_{k}^+)$.

 Our next objective is to show that the system of intervals
 $\{(a_k,b_k)\}$ constructed above satisfies (\ref{1.6}) and
 (\ref{1.7}).

 By (\ref{4.7})-(\ref{4.8}) with $a=a_k$ and $b=b_k$ there exists $0<q<1$ such that
 $$\frac{\kap(E\cap[a_k,b_k])}{\kap([a_k,b_k])}<q,
 $$
 from which (\ref{1.6}) follows.

 Let $\Ga_k$ be the family of all crosscuts of $\He$ which
 separate points $a_k$ and $b_k$ from $0$ and $\infty$. By (\ref{3.nn4}) for its
 module we have
 \begin{eqnarray}
 m(\Ga_k)&\ge& m(\{\ga_r:=\{z\in\He:\,|z-b_k|=r\}:\,
 b_k-a_k<r<b_k\})\nonumber\\
 \label{6.c1}&=&\frac{1}{\pi}\log\frac{b_k}{b_k-a_k}.v_{j_k}
 \end{eqnarray}
 Let $\Ga_k'=\phi(\Ga_k)$. Consider the metric
 $$\rho(w):=\left\{\begin{array}{ll}
 1,&\mb{ if }-v_{j_k}\le\Re w\le u_{j_k}+v_{j_k},0\le\Im w\le 2v_{j_k},\\[2ex]
 0,&\mb{ elsewhere}.
\end{array}\right.$$
Since
 $$\int_\ga\rho(w)|dw|\ge v_{j_k}\quad (\ga\in\Ga'_k),
 $$
 it follows that
 \beq\label{6.c2}
 m(\Ga_k')\le\frac{A(\rho)}{v^2_{j_k}}=\frac{2(2v_{j_k}+u_{j_k})}{v_{j_k}}
 \le2\left(2+\frac{1}{c}\right)=:C.
 \eeq
 Comparing (\ref{6.c1}) and (\ref{6.c2}) we obtain
 $$\frac{b_k-a_k}{b_k}\ge e^{-\pi C},
 $$
 from which
 (\ref{1.7}) follows.

 \hfill\rbox

\absatz{Other Proofs}

{\it Proof of Remark 3.} Our objective is to construct a cover of
$E^*$ satisfying the conditions of Remark 1.

First, we divide $\R$ by points $\pm 2^k,k=0,1,\ldots$ on the
infinite number of intervals of the form either
$I_k^+=[2^k,2^{k+1}]$, or $I_k^-=[-2^{k+1},-2^{k}]$, or
$I_0=I=[-1,1]$ (here $k=1,2,\ldots).$

For each interval $I'=I_k^+$ with the property
 \beq\label{15.1}
 |I'\cap E^*|\ge 2^{k-1},
 \eeq
 we have
 $$
 \int_{2^k}^{2^{k+1}}\frac{\theta^2_E(t)}{t^3}dt\ge
 \int_{3\cdot2^{k-1}}^{2^{k+1}}\frac{(t-3\cdot2^{k-1})^2}{t^3}dt
 =\log\frac{4}{3}-\frac{9}{32}>10^{-3}.
 $$
 The same inequality is valid for $I'$ of the form $I'=I_k^-$
 satisfying (\ref{15.1}). Therefore, there is only a finite number of
 intervals, satisfying (\ref{15.1}). According to Theorem C
 we can assume that $E^*\subset\cup_{k=k_0}^\infty I_k^\pm$ for
 sufficiently large $k_0$ such that
 $$|I_k^\pm\cap E^*|<\frac{1}{2}|I_k^\pm|\quad (k\ge k_0).
 $$
 For any interval of the covering of $E^*$ (say $I_k^+$) which has nonempty
 intersection with $E^*$,
 the linear transformation of
the result of the lemma in \cite[p. 580]{and} (with its obvious
extension to the case of sets consisting of the infinite number of
intervals) implies the existence of real numbers
 $$2^k\le c_{k,1}^+<d_{k,1}^+\le c_{k,2}^+<\ldots\le c_{k,n_k^+}^+<d_{k,n_k^+}^+\le 2^{k+1}
 $$
  such that
 $$I_k^+\cap E^*\subset \bigcup_{j=1}^{n_k^+}[c_{k,j}^+,d_{k,j}^+],
 $$
  $$ |[c_{k,j}^+,d_{k,j}^+]\cap
 E^*|\ge \frac{e_{k,j}^+}{4},\quad e_{k,j}^+:=d_{k,j}^+-c_{k,j}^+.$$
 Notice that
 \begin{eqnarray*}
\sum_{j=1}^{n_k^+}\left(\frac{d_{k,j}^+-c_{k,j}^+}{|c_{k,j}^+|+1}\right)^2
 &\le&2^{-2k} \sum_{j=1}^{n_k^+}e_{k,j}^{+2}\le
2^{-2k}\left( \sum_{j=1}^{n_k^+}e_{k,j}^+\right)^2\\
 &\le& 2^{4-2k}\theta^2_E(2^{k+1})\le
 2^8\int\limits_{2^{k+1}}^{2^{k+2} }\frac{\theta^2_E(t)}{t^3}dt\,
 .\end{eqnarray*}
 Consider the system of intervals $\{\{(c_{k,j}^\pm,d_{k,j}^\pm)\}_{j=1}^{n_k^\pm}\}_{k\ge
 k_0}$ constructed as above for each interval $I_k^\pm,k\ge k_0$ which has nonempty intersection
 with $E^*$.
 To use Remark 1 we need to cover $E^*$ by an infinite number of
 intervals going in both positive and negative directions of $\R$.
 To satisfy this formal condition we can add to the intervals
 constructed above the intervals of the form $(n,n+1)\subset E$ or
 $(-n-1,-n)\subset E$, where $n$ is a positive integer if there is
 only a finite number of $I_k^\pm$ such that $I_k^\pm\cap
 E^*\neq\emptyset$.

 As a result we obtain the cover of $E^*$ satisfying
 (\ref{1.41n}), (\ref{1.42n}), (\ref{1.5}), and (\ref{1.8}).

 \hfill\rbox

 {\it Proof of Remark 4.}
 Following the idea from \cite[Section 3.1]{tot} we set
 $E=\R\setminus E^*$, where
 $$E^*=\bigcup_{j=1}^\infty\left(2^j,2^j+\frac{1}{4}\left(\theta(2^{j})-
 \theta(2^{j-1})\right)\right)=: \bigcup_{j=1}^\infty(c_j,d_j).$$
 Since for $j\ge 1$ and $2^j<t\le 2^{j+1},$
 $$\theta_E(t)\le\sum_{k=1}^j \frac{1}{4}\left(\theta(2^{k})-
 \theta(2^{k-1})\right)=\frac{1}{4}\theta(2^j)\le\theta(t),$$
 we obtain (\ref{1.3n}).

 Setting $s_j:=2^{-j}\theta(2^j)$, for any $n>1$ we have
 \begin{eqnarray}
 \sum_{j=1}^n\left(\frac{d_j-c_j}{|c_j|+1}\right)^2&\ge&
 2^{-6}\sum_{j=1}^n\left(\frac{\theta(2^j)-\theta(2^{j-1})}{2^j}\right)^2
 =
 2^{-6}\sum_{j=1}^n\left(s_j-\frac{1}{2}s_{j-1}\right)^2\nonumber\\
 &=&2^{-6}\left(\sum_{j=1}^ns_j^2- \sum_{j=1}^ns_js_{j-1}+\frac{1}{4}
 \sum_{j=1}^ns_{j-1}^2\right)\nonumber\\
 \label{1.4n}
 &\ge&
 2^{-8}\sum_{j=1}^ns_{j-1}^2-2^{-4},
 \end{eqnarray}
 where in the last step we used the estimate
 $$\sum_{j=1}^ns_j^2- \sum_{j=2}^{n}s_js_{j-1}\ge 0,$$
 following from the Cauchy-Schwarz inequality, and the estimate
 $s_1s_{0}\le 4$, following from the assumption (\ref{1.hu1}).
 Furthermore,
 \beq\label{1.5n}
 \int_1^{2^{n+1}}\frac{\theta^2(t)}{t^3}dt=\sum_{j=0}^n\int_{2^j}^{2^{j+1}}
 \frac{\theta^2(t)}{t^3}dt\le
 \sum_{j=0}^n\frac{\theta^2(2^{j+1})}{2^{2j}}=4\sum_{j=0}^n
 s^2_{j+1}.
 \eeq
 Thus, (\ref{1.4n}), (\ref{1.5n}) and (\ref{1.hu2})  imply (\ref{1.10}), i.e.,
 $$\sum_{j=1}^\infty\left(\frac{d_j-c_j}{|c_j|+1}\right)^2=\infty.
 $$
  To complete the proof  we  apply
  Remark 2.

 \hfill\rbox

{\it Acknowledgement.}
 The
author would like to warmly thank  M. Nesterenko for many useful
remarks.

Vladimir V. Andrievskii

 Department of Mathematical Sciences

 Kent State University

 Kent, OH 44242, U.S.A.

andriyev@math.kent.edu

\end{document}